\documentclass[11pt]{amsart}
\usepackage[T1]{fontenc}
\usepackage[applemac]{inputenc}
\usepackage{amssymb}
\usepackage{amsmath}
\usepackage{amsthm}
\usepackage{amsfonts}
\usepackage[all]{xy}
\usepackage{lscape}
\usepackage{vmargin}
    \setpapersize{A4}
    \setmargins{2.9cm}{2.5cm}{15.2cm}{23.2cm}{12pt}{25pt}{0pt}{30pt}

\theoremstyle{definition}
\newtheorem*{rem}{Remark}
\newtheorem{prop}{Proposition}
\newtheorem*{fact}{Fact}
\newtheorem{thm}{Theorem}
\newtheorem*{defn}{Definition}

\newcommand{\z}{\mathbb{Z}}
\newcommand{\q}{\mathbb{Q}}
\newcommand{\rea}{\mathbb{R}}
\newcommand{\co}{\mathbb{C}}
\newcommand{\rw}{\rightarrow}
\newcommand{\lrw}{\longrightarrow}
\newcommand{\ad}{\mathbb{A}}
\newcommand{\af}{\mathbb{A}_{f}}
\newcommand{\A}{\mathbf{A}}

\newcommand{\G}{\mathbf{G}}
\newcommand{\qH}{\mathbf{H}}
\newcommand{\qL}{\mathbf{L}}
\newcommand{\M}{\mathbf{M}}
\newcommand{\N}{\mathbf{N}}
\newcommand{\qP}{\mathbf{P}}
\newcommand{\qQ}{\mathbf{Q}}
\newcommand{\tS}{\mathbf{S}}
\newcommand{\T}{\mathbf{T}}

\makeatother

\begin{document}
\title{Eisenstein cohomology for congruence subgroups of $SO(n,2)$}
\author{Gerald Gotsbacher}
\address{Department of Mathematics, University of Toronto, Toronto, ON, M5S 2E4
Canada}
\email{gerald@math.toronto.edu}
\thanks{The author was supported by the Austrian Science Fund (FWF), project no. P16762 at the University of Vienna when the results of this paper were obtained.}

\begin{abstract}
The automorphic cohomology of a connected reductive algebraic group defined over $\q$
decomposes as a direct algebraic sum of cuspidal and Eisenstein cohomology.  
The present paper investigates the Eisenstein cohomology for congruence subgroups of a rational form $\G$ of $\q$-rank $2$ of $SO(n,2)$ in the generic case. The main result provides a description of the internal structure of the summands in the Eisenstein cohomology corresponding to maximal parabolic $\q$-subgroups. 
\end{abstract}
\maketitle

\tableofcontents

\section*{Introduction}

Let there be given a connected linear semisimple algebraic group $\G$ defined over $\q$, a congruence subgroup $\Gamma\subset\G(\q)$ and an irreducible finite dimensional complex representation $(\tau,E)$ of $\Gamma$. 
The group cohomology $H^*(\Gamma,E)$ of $\Gamma$ with coefficients in $E$ has a description
in terms of the automorphic spectrum of $\Gamma$.
More precisely, let $\mathfrak{g}$ denote the Lie algebra of $\G(\co)$, $K$ some maximal compact subgroup of $G=\G(\rea)$ and $\mathcal{A}_E(\Gamma\backslash G)$ the space of automorphic forms on $\Gamma\backslash G$ with respect to $E$. It is then a consequence of a result of J. Franke (Theorem $18$ in \cite{Fra}) that the group cohomology $H^*(\Gamma,E)$ is isomorphic to the automorphic cohomology $H^*(\mathfrak{g},K; \mathcal{A}_E(\Gamma\backslash G)\otimes E)$ of $\Gamma$.

\textit{Eisenstein cohomology} refers to the attempt to construct cohomology classes (in the automorphic cohomology of $\Gamma$) in terms of Eisenstein series associated to cusp forms for the projection of 
$\Gamma$ to a Levi subgroup of $\G$ and to ultimately describe
a natural orthogonal complement of the cuspidal cohomology in the automorphic cohomology of 
$\Gamma$. The basic idea is due to G. Harder (cf. \cite{Har2}) and appeals to the theory of Eisenstein series as developed by R. Langlands (cf. \cite{Lan}, \cite{MW}). 

The remainder of this introduction gives a brief outline of the technique of Eisenstein cohomology in terms of the algebraic $\q$-group $\G$ chosen for the present paper, though in the ad\`elic and representation-theoretic setting.\footnote{For a thorough 
survey of the field see \cite{Sch3}, for a detailed description of the method see \cite{Sch1}.}

Thus, for the moment $\G$ is to denote a $\q$-rational form of $\q$-rank $2$ of the real semisimple Lie group $SO(n,2)$ for $n\geq5$, $(\tau,E)$ to be an irreducible finite-dimensional complex representation
of $\G$ and $\mathcal{A}_{E}=\mathcal{A}_{E}(\G(\q)\backslash\G(\ad))$ the space of ad\`elic automorphic forms for $\G$ with respect to $E$.
The automorphic cohomology of $\G$ with coefficients in $E$ is the relative Lie algebra cohomology
$H^*(\mathfrak{g},K;\mathcal{A}_{E}\otimes E)$ allowing of a decomposition as an algebraic direct sum
$$H^*(\mathfrak{g},K;{}^{\circ}\mathcal{A}_{E}\otimes E)\oplus
\bigoplus_{\{\qP\}\neq\{\G\}}H^*(\mathfrak{g},K;\mathcal{A}_{E,\{\qP\}}\otimes E),$$
where ${}^{\circ}\mathcal{A}_{E}$ denotes the space of cusp forms in $\mathcal{A}_{E}$, $\{\qP\}$ the 
associate class of a parabolic $\q$-subgroup $\qP\subset\G$ and $\mathcal{A}_{E,\{\qP\}}\subset\mathcal{A}_{E}$ the space of automorphic forms negligible along every parabolic $\q$-subgroup $\qQ\subset\G$ with $\qQ\not\in\{\qP\}$. The algebraic direct sum of the summands indexed by $\{\qP\}\neq\{\G\}$ is called the \textit{Eisenstein cohomology} of $\G$.

As for the scope of the present work the account of the Eisenstein cohomology for maximal parabolic $\q$-subgroups given is subject to the overall hypothesis that the highest weight $\lambda_{\tau}$ of $(\tau,E)$ be regular. That is, it is \textit{regular} Eisenstein cohomology classes which are constructed.
Since the associate and conjugacy classes of parabolic $\q$-subgroups coincide for the given group $\G$, it suffices to restrict to standard parabolic $\q$-subgroups $\qP$. 

The construction starts from so-called classes of type $(\pi,w)$ in the cuspidal cohomology of the Levi component $\qL_{\qP}$ of  $\qP$. 
The type comprises a minimal coset representative $w$ for the orbit space of the Weyl group $W_{\qL_{\qP}}$ of $\qL_{\qP}$ acting on the Weyl group $W$ of $\G$; and essentially an irreducible unitary representation $\pi$ of the group of real points of $\qL_{\qP}$ with non-trivial cohomology twisted
by an irreducible finite-dimensional complex rational representation of $\qL_{\qP}$ induced by $(\tau,E)$ and depending on $w$. 
The set $W^{\qP}$ of minimal coset representatives $w$ is determined explicitly for the two standard maximal parabolic $\q$-subgroups $\qP\subset\G$ by way of producing what is called the \emph{Hasse diagram} of a partially ordered set, in this case of the $W$-orbit of a certain weight associated to $\qP$. 
Section \ref{Kostant} gives the results. 
In particular, it is the usage of the Hasse diagram of $\qP$ enabling one to overcome the technical difficulties in dealing with the elements of $W^{\qP}$ arising from the fact that its order grows 
polynomially along with $n$. Appendix \ref{diagram} provides a concise introduction to the notion of Hasse diagram as it is utilised here.
The representations $\pi$ of the Levi components $\qL_{\qP}$ for the two standard maximal parabolic $\q$-subgroups $\qP\subset\G$ are described in section \ref{cohrep} up to unitary equivalence according to the Vogan-Zuckerman classification. 

Finally, lifting a given class of type $(\pi,w)$ into $H^*(\mathfrak{g},K;\mathcal{A}_{E,\{\qP\}}\otimes E)$ by means of parabolic induction followed by Eisenstein summation, the thus obtained Eisenstein series is evaluated at a certain point $\lambda_w$ stemming from $\lambda_{\tau}$ in a way that is uniquely determined by $w\in W^{\qP}$ and all non-trivial classes are identified. 
A description of the summands $H^*(\mathfrak{g},K;\mathcal{A}_{E,\{\qP\}}\otimes E)$ in the Eisenstein cohomology corresponding to standard maximal parabolic $\q$-subgroups $\qP\subset\G$ is displayed in the main Theorems \ref{myfirst} to \ref{mythird} of the paper. 

I'd like to express my gratitude to Joachim Schwermer for sharing his knowledge and expertise.

\section{Preliminaries}\label{pre}

\subsection*{Reductive algebraic $\q$-groups}
As in the introduction $\G$ will always denote some connected linear algebraic group defined over $\q$, which henceforth is allowed to be reductive.
For any commutative $\q$-algebra $A$ the group of $A$-points of $\G$ is defined to be $\G(A)=\G\cap\mathbf{GL}_n(A)$. 

The group of complex, resp. rational characters of $\G$ defined over $\q$ is denoted as $X(\G)$, resp. $X_{\q}(\G)$. The algebraic $\q$-subgroup ${}^{\circ}\G=\bigcap_{\chi\in X_{\q}(\G)}\textrm{ker}\chi^2$ is normal in  $\G$ and contains all compact and arithmetic subgroups of $\G(\rea)$.

\subsubsection*{Absolute root system and minimal coset representatives} 
The root system for a pair $(\G,\T)$ where $\T$ is a maximal $\q$-torus in $\G$, is denoted as
$\Phi=\Phi(\G,\T)$. The Weyl group $W(\G,\T)=\N_{\G}(\T)/\mathbf{Z}_{\G}(\T)$ of the pair $(\G,\T)$ is frequently identified to the Weyl group $W=W(\Phi)$ of the root system $\Phi$ without mention.

Let $\mathbf{B}\subset\G$ be a Borel subgroup such that $\mathbf{B}\supset\T$ with corresponding  positive system $\Phi^+$. 
Suppose $\qP\supset\mathbf{B}$ be a standard parabolic subgroup of $\G$ and $\qL_{\qP}\subset\qP$ some Levi subgroup. The 
Weyl group $W_{\qL_{\qP}}=W(\Phi(\qL_{\qP},\T))$ of the root system $\Phi(\qL_{\qP},\T)$ can be embedded into 
$W$ and the quotient $W_{\qL_{\qP}}\backslash W$ has a distinguished set of representatives, the set
$$W^{\mathbf{P}}=\{w\in W\vert w^{-1}(\Delta(\qL_{\qP},\T))\subset\Phi^+\}$$
of \textit{minimal coset representatives}. Here, $\Delta(\qL_{\qP},\T)$ is to denote the set of simple roots of the positive system on $\Phi(\qL_{\qP},\T)$
induced by $\Phi^+$. 

\subsubsection*{Relative root system and parabolic $\q$-subgroups}
Given a \textit{maximal $\q$-split} torus in $\tS\subset\G$ such that  $\tS\subset\T$
the set of rational roots is denoted as ${}_{\q}\Phi={}_{\q}\Phi(\G,\tS)$. 
Moreover, if ${}_{\q}\Phi^+$ is a positive system it is assumed that $\Phi^+$ be a compatible positive system, that is ${}_{\q}\Phi^+$ lies in the image of $\Phi^+$ under the restriction map $X(\T)\rw X(\tS)$. 

The parabolic $\q$-subgroups $\qP$ which are standard with respect to $\tS$ and ${}_{\q}\Phi^+$ are in
$1-1$ correspondence with the subsets $I$ of the set ${}_{\q}\Delta\subset{}_{\q}\Phi^+$ of simple rational roots.
Specifically, if $I\subset{}_{\q}\Delta$, set $\A_{I}=(\bigcap_{\alpha\in I}\ker\alpha)^{\circ}$ for the subtorus of $\tS$ defined by $I$ and $\qL_{I}=\mathbf{Z}_{\G}(\A_{I})$ for its centraliser. The latter one is the Levi component of a unique parabolic $\q$-subgroup $\qP_{I}$ of $\G$ with unipotent radical denoted by $\N_{I}$. Since $\qL_I$ is an almost direct product $\qL_{I}=\M_{I}\A_{I}$ where $\M_{I}={}^{\circ}\qL_{I}$, the Levi decomposition of $\qP_I$ admits the refinement $\qP_{I}=\M_{I}\A_{I}\N_{I}$. 
For an arbitrary parabolic $\q$-subgroup $\qP$ of $\G$ the alike decomposition is written as 
$\qP=\M_{\qP}\A_{\qP}\N_{\qP}$.

Two parabolic $\q$-subgroups $\qP, \qQ\subset\G$
are called \textit{associate}, if there is an element $g\in\G(\q)$ such that
${}^g\qL_{\qQ}=\qL_{\qP}$, where ${}^g{\cdot}=\textrm{int}(g)$ denotes the inner automorphism on $\G$ induced by $g$. The associate class of $\qP$ is denoted by $\{\qP\}$ and the finite collection of all these by ${}\mathcal{C}$.

\subsection*{The group of real points and Lie algebras}
The complex Lie algebra of an algebraic $\q$-group $\G$ will be denoted as $\mathfrak{g}=\textrm{Lie}(\G(\co))$. The real form of  $\mathfrak{g}$ corresponding to the group $\G(\rea)$ of real points will be indicated by putting a `$0$' as subscript, that is $(\mathfrak{g}_0)_{\co}=\mathfrak{g}$. In addition, for an abelian Lie algebra $\mathfrak{a}$ the dual is denoted as $\mathfrak{\check{a}}$. If $\mathfrak{h}\subset\mathfrak{g}$ is a Cartan subalgebra an arbitrary element in $\mathfrak{\check{h}}$ will be denoted by $\lambda$, as will be its restriction to $\mathfrak{a}_{\qP}\subset\mathfrak{h}$.

Often the group $\G(\rea)$ of real points of an algebraic $\q$-group will be denoted by $G$, the connected component containing the identity by $G^+$. The symmetric space $X$ associated to 
$G$ is the space of maximal compact subgroups of $G$.
If $K\subset G$ is a maximal compact subgroup let $\Theta$, resp. $\theta$ denote the corresponding Cartan involution on $G$, resp. $\mathfrak{g}_0$. The Cartan decomposition is written as $\mathfrak{g}_0=\mathfrak{k}_0\oplus
\mathfrak{p}_0$ with $\mathfrak{k}_0=\textrm{Lie}K$ and $\mathfrak{p}_0$ the $-1$-eigenspace of $\theta$. Whereas the intersection of $K$ with a subgroup $H\subset G$ will be written as $K_H$, the 
corresponding restrictions of the Cartan involution, resp. decomposition are denoted with the same letters. 
By common abuse of language the unique complex linear extension of $\theta$ to $\mathfrak{g}$
is called Cartan involution again and denoted by the same letter. In addition, let $\sigma$ denote the conjugation on $\mathfrak{g}$ induced by $\mathfrak{g}_0$.

Let $\mathfrak{h}_0\subset\mathfrak{g}_0$ be a maximally compact $\theta$-stable Cartan
subalgebra and set $\mathfrak{t}_0=\mathfrak{h}_0\cap\mathfrak{k}_0$ for the compact, 
$\mathfrak{a}_0=\mathfrak{h}_0\cap\mathfrak{p}_0$ for the non-compact part. 
Next, suppose $\mathfrak{q}\subset\mathfrak{g}$ is a \emph{$\theta$-stable }parabolic subalgebra
containing $\mathfrak{h}$, i.e. $\mathfrak{q}$ is stable under $\theta$ and $\sigma\mathfrak{q}\cap\mathfrak{q}$ is a maximal reductive Lie subalgebra. Let $\mathfrak{u}$ be its nil radical and $\mathfrak{q}=\mathfrak{l}\oplus\mathfrak{u}$ a Levi decomposition.  In \ref{cohrep} the set $\Phi(\mathfrak{u}\cap\mathfrak{p})$ of weights for the adjoint action of $\mathfrak{t}$ in $\mathfrak{u}\cap\mathfrak{p}$ will be of significance.
In addition, $\rho(\mathfrak{u}\cap\mathfrak{p})$ is to denote half the sum of all weights in $\Phi(\mathfrak{u}\cap\mathfrak{p})$. Notation will be analogous for any reductive Lie subalgebra of $\mathfrak{g}_0$.

\subsection*{The group of ad\`eles and ad\`elic automorphic forms}
The ring of ad\`eles over $\q$ is denoted as $\ad$. It is the direct product $\ad=\rea\times\af$ of the
field of real numbers and the ring $\af$ of finite ad\`eles.

Let $K_{\mathrm{ad}}\subset\G(\ad)$ be a maximal compact subgroup of the group of ad\`eles of $\G$. 
It writes as a product $K_{\mathrm{ad}}=KK_f$ with $K_f\subset\G(\af)$ a restricted product $K_f=\prod_{p}K_p$ over all primes $p\in\q$, and $K_p\subset\G(\q_p)$ a maximal compact subgroup. In addition, it is assumed that $K_{\mathrm{ad}}$ is in \emph{good position} relative to the minimal parabolic $\q$-subgroup 
$\qP_0\subset\G$ corresponding to the choice of ${}_{\q}\Phi^+$. The latter condition implies that the standard height function $H_{\qP}:\qP(\ad)\rightarrow\mathfrak{a}_{\qP}$ for a standard parabolic $\q$-subgroup $\qP\supset\qP_0$ has an extension to $\G(\ad)$, again denoted by $H_{\qP}$.

Let $\mathcal{Z}(\mathfrak{g})$ be the center of the universal enveloping algebra $\mathcal{U}(\mathfrak{g})$ of $\mathfrak{g}$. Then $\mathcal{Z}(\mathfrak{g})$ acts on the dual representation $\check{E}$ of $E$.
Let $\mathcal{I}\subset\mathcal{Z}(\mathfrak{g})$ be the annihilator of $\check{E}$ in $\mathcal{Z}(\mathfrak{g})$ and set $\mathcal{A}_E$ for the space of complex-valued, $K_{\mathrm{ad}}$-finite, smooth functions of uniform moderate growth on $\G(\q)\backslash\G(\ad)$ which are annihilated by a power of $\mathcal{I}$.
The space $\mathcal{A}_E$ is a $(\mathfrak{g},K;\G(\af))$-module and its elements are called ad\`elic automorphic forms for $\G$ (relative to $E$).

\section{Automorphic cohomology for congruence groups}\label{automorphic}

In the present and subsequent section $\G$ is assumed to have $\q$-anisotropic center and positive semisimple $\q$-rank. As before $(\tau,E)$ denotes some irreducible finite-dimensional complex representation of $\G$ of highest weight $\lambda_{\tau}$.
As stated in the introduction the automorphic cohomology of $\G$ with respect to $E$ by definition is the relative Lie algebra cohomology $H^*(\mathfrak{g},K;\mathcal{A}_E\otimes E)$ with coefficients in $\mathcal{A}_E$ twisted by $E$.

\subsection{Decomposition along the cuspidal support}
 The space $\mathcal{A}_E$ allows of a decomposition along the cuspidal support of distinct associate classes of parabolic $\q$-subgroups. (See \cite{F-S}, \cite{MW} for details).

First, let $V_{\G}=C^{\infty}_{\textrm{umg}}(\G(\q)\backslash\G(\ad))$ be the space of smooth, complex-valued functions on $\G(\q)\backslash\G(\ad)$ of \textit{uniform moderate growth}, and for $\{\qP\}\in\mathcal{C}$ set $V_{\G}(\{\qP\})$ for the space of elements of $V_{\G}$ which are \textit{negligible}
along $\qQ\subset\G$ for every parabolic $\q$-subgroup $\qQ\not\in\{\qP\}$. That is, the constant term of an element in $V_{\G}(\{\qP\})$ with respect to $\qQ$ is orthogonal to the space of cusp forms on $\M_{\qQ}$. The space $V_{\G}$ is a $(\mathfrak{g},K;\G(\af))$-module and admits a decomposition $$V_{\G}=\bigoplus_{\{\qP\}\in\mathcal{C}}V_{\G}(\{\qP\})$$ as a finite direct sum of $(\mathfrak{g},K;\G(\af))$-modules as proved by Langlands\footnote{See A. Borel, J.-P. Labesse, J. Schwermer, \emph{On the cuspidal cohomology of $S$-arithmetic subgroups of reductive groups over number fields}.  Compositio Math.  \textbf{102}  (1996),  no. 1, 1--40.}.
This decomposition descends to the submodule $\mathcal{A}_E$, that is $$\mathcal{A}_E=\bigoplus_{\{\qP\}\in\mathcal{C}}\mathcal{A}_{E,\{\qP\}}$$ with $\mathcal{A}_{E,\{\qP\}}=\mathcal{A}_E\cap V_{\G}(\{\qP\})$ and the isomorphism being one of $(\mathfrak{g},K;\G(\af))$-modules. (Of course, 
$\mathcal{A}_{E,\{\G\}}={}^{\circ}\mathcal{A}_E$).

Next, let $\qQ\in\{\qP\}$ and $\pi$ be an irreducible representation of $\qL_{\qQ}(\ad)$, unitary modulo the center, and such that $(1)$ the central character $\chi_{\pi}:\A_{\qQ}(\ad)\rw\co^{\times}$ is trivial on $\A_{\qQ}(\q)$, $(2)$ $\pi$ occurs in the cuspidal summand ${}^{\circ}L^2(\qL_{\qQ}(\q)\backslash\qL_{\qQ}(\ad))_{\chi_{\pi}}$ and $(3)$ the infinitesimal character of $\pi$ matches the infinitesimal character of the dual representation $\check{E}$. These conditions entail some compatibility requirements (cf. 
\cite{F-S} p. 771) and we let $\phi_{\qQ}\ni\pi$ be a finite set of such representations meeting them. As a
consequence, two such sets $\phi_{\qQ},\phi_{\qQ'}$ for $\qQ,\qQ'\in\{\qP\}$ are associate by the dual of the inner automorphism of $\G$ mapping the respective Levi components to one another. Finally, 
let $\Phi_{E,\{\qP\}}$ denote the collection of all classes $\phi=\{\phi_{\qQ}\}_{\qQ\in\{\qP\}}$.

Let $\pi\in\phi_{\qQ}$ for some $\phi\in\Phi_{E,\{\qP\}}$, where $\qQ\in\{\qP\}$ and set $d\chi_{\pi}\in\mathfrak{\check{a}}_{\qQ}$ for the differential of the central character $\chi_{\pi}$ restricted to $A_{\qQ}^+$. Then ${}^{\circ}L^2_{\pi}(\qL_{\qQ}(\q)A_{\qQ}^+\backslash\qL_{\qQ}(\ad))_{\chi_{\pi}}$ shall denote the space of all cuspidal automorphic forms on $\qL_{\qQ}(\ad)$ which transform according to $\pi$. Furthermore, let $\tilde\pi$ denote the irreducible unitary representation of $\qL_{\qQ}(\ad)$ 
obtained by normalising $\pi$ with respect to the action of the central character $\chi_{\pi}$ on 
$A_{\qQ}$. Let $W_{\qQ,\tilde\pi}$ be the $(\mathfrak{g},K)$-module of smooth, $K_{\mathrm{ad}}$-finite functions $f:\qL_{\qQ}(\q)\N_{\qQ}(\ad)A_{\qQ}^+\backslash\G(\ad)\rw\co$ such that for any $g\in\G(\ad)$ the function $l\mapsto f(lg)$ for $l\in\qL_{\qQ}(\ad)$ belongs to 
${}^{\circ}L^2_{\tilde\pi}(\qL_{\qQ}(\q)A_{\qQ}^+\backslash\qL_{\qQ}(\ad))_{\chi_{\pi}}$.
To $f\in W_{\qQ,\tilde\pi}$ and $\lambda\in\mathfrak{\check{a}}_{\qQ}$ we now associate an Eisenstein
series $E^{\G}_{\qQ}(f,\lambda)$ to be a function in $g\in\G(\ad)$, whenever convergent:
$$E^{\G}_{\qQ}(f,\lambda)(g)=\sum_{\gamma\in\qQ(\q)\backslash\G(\q)}e^{\langle H_{\qQ}(\gamma g),\lambda+\rho_{\qQ}\rangle}f(\gamma g).$$
Given that the real part of $\lambda$ lies inside the positive Weyl chamber defined by 
$\qQ$ being in addition sufficiently regular this series is known to converge normally for $g$ in a compact set. As a function in $\lambda$ it admits a meromorphic continuation to all of $\mathfrak{\check{a}}_{\qQ}$ (cf. \cite{MW} II.1.5).

Finally, we define $\mathcal{A}_{E,\{\qP\},\phi}\subset\mathcal{A}_{E,\{\qP\}}$ as the space generated by all residues and derivatives of Eisenstein series $E^{\G}_{\qQ}(f,\lambda)$ for $f$ ranging through $W_{\qQ,\tilde\pi}$.
(Notice that this space is denoted by $\tilde{\mathcal{A}}_{E,\{\qP\},\phi}$ in \cite{F-S}).
Then there is a decomposition 
$$\mathcal{A}_E=\bigoplus_{\{\qP\}\in\mathcal{C}}\bigoplus_{\phi\in\Phi_{E,\{\qP\}}}\mathcal{A}_{E,\{\qP\},\phi}$$ of $(\mathfrak{g},K;\G(\af))$-modules giving rise to the 

\begin{thm}(\cite{F-S})
Let $\G$ be a connected reductive algebraic $\q$-group the center of which is anisotropic over $\q$, and 
suppose it have positive semisimple $\q$-rank. Then there is a direct sum decomposition in cohomology
$$H^*(\mathfrak{g},K;\mathcal{A}_{E}\otimes E)=\bigoplus_{\{\qP\}\in\mathcal{C}}\bigoplus_{\phi\in\Phi_{E,\{\qP\}}} H^*(\mathfrak{g},K;\mathcal{A}_{E,\{\qP\},\phi}\otimes E),$$
where as before $\mathcal{C}$ denotes the set of classes $\{\qP\}$ of associate parabolic $\q$-subgroups and $\Phi_{E,\{\qP\}}$ the set of classes $\phi=\{\phi_{\qQ}\}_{\qQ\in\{\qP\}}$ of associate irreducible cuspidal automorphic representations of the Levi components of elements $\qQ\in\{\qP\}$. 
\end{thm}

\subsection{Cuspidal classes of type $(\pi,w)$}\label{type}
Notation is as before, except henceforth we choose $\qP\in\{\qP\}$ and $\pi$ shall denote either a cuspidal automorphic representation of $\qL_{\qP}$ belonging to some $\phi\in\Phi_{E,\qP}$ or its archimedean component. For $\lambda\in\mathfrak{\check{a}}_{\qP}$ let $\textrm{Ind}^{\G(\ad)}_{\qP(\ad),\pi,\lambda}$ denote the representation of $\G$ induced from $\pi$ by parabolic induction. Then there is an isomorphism $$W_{\qP,\tilde\pi}\otimes\lambda\simeq (\textrm{Ind}^{\G(\ad)}_{\qP(\ad),\pi,\lambda})^{m_0(\pi)}$$
of $(\mathfrak{g},K)$-modules, where $m_0(\pi)$ denotes the multiplicity of $\pi$ in ${}^{\circ}L^2(\qL_{\qP}(\q)\backslash\qL_{\qP}(\ad))$. This way the $\G(\af)$-action on $\textrm{Ind}^{\G(\ad)}_{\qP(\ad),\pi,\lambda}$ carries over to $W_{\qP,\tilde\pi}$. Moreover, the symmetric algebra $\mathrm{Sym}(\mathfrak{\check{a}}_{\qP})$ of $\mathfrak{\check{a}}_{\qP}$ can be regarded as the space of differential operators $\frac{\partial^{\nu}}{\partial\lambda^{\nu}}$ with constant coefficients in $\mathfrak{\check{a}}_{\qP}$ for some multi-index $\nu$.  

Now, there exists a polynomial function $q(\lambda)$ on $\mathfrak{\check{a}}_{\qP}$ such that for every $f\in W_{\qP,\tilde\pi}$ there is a neighbourhood of $d\chi_{\pi}\in\mathfrak{\check{a}}_{\qP}$,
in which $q(\lambda)E^{\G}_{\qP}(f,\lambda)$ is holomorphic (cf. \cite{MW} IV.1). 
Hence, for $\pi\in\phi_{\qP}$ and $\phi_{\qP}\in\phi$ the mapping 
$$f\otimes\frac{\partial^{\nu}}{\partial\lambda^{\nu}}\longmapsto\frac{\partial^{\nu}}{\partial\lambda^{\nu}}
(q(\lambda)E^{\G}_{\qP}(f,\lambda))_{|d\chi_{\pi}}$$
which assigns to an element $f\otimes\frac{\partial^{\nu}}{\partial\lambda^{\nu}}$  
the derivative of $q(\lambda)E^{\G}_{\qP}(f,\lambda)$ at $d\chi_{\pi}$ with respect to 
$\frac{\partial^{\nu}}{\partial\lambda^{\nu}}$ 
yields a homomorphism $$W_{\qP,\tilde\pi}\otimes\mathrm{Sym}(\mathfrak{\check{a}}_{\qP})\rw\mathcal{A}_{E,\{\qP\},\phi}$$
of $(\mathfrak{g},K;\G(\af))$-modules. 

The Hochschild-Serre spectral sequence in the category of $(\mathfrak{p},K_P)$-modules\footnote{This is the only instance where $\mathfrak{p}=\mathrm{Lie}(\qP(\co))$.} (cf. \cite{B-W} III. Thm 3.3) associated to the double complex obtained from replacing the coefficient module in
$D^*=D^*(\mathfrak{p},K_{P}; H_{\pi}\otimes\mathrm{Sym}(\mathfrak{\check{a}}_{\qP})\otimes E)$ by some injective resolution of it abuts to the cohomology of $D^*$, degenerates and has second term $E_2^{p,q}=H^p(\mathfrak{l}_{\qP},K_{L_{\qP}};H_{\pi}\otimes H^q(\mathfrak{n}_{\qP},E)\otimes\mathrm{Sym}(\mathfrak{\check{a}}_{\qP}))$.
Provided that the highest weight $\lambda_{\tau}$ of $(\tau,E)$ is dominant by a theorem of Kostant (cf. \cite{Kos}, 5.13) there is an isomorphism of $\qL_{\qP}(\co)$-modules
$$H^q(\mathfrak{n}_{\qP},E)=\bigoplus_{w\in W^{\qP}, \,l(w)=q}F_{\mu_w}$$
where $\mu_w=w(\lambda_{\tau}+\rho)-\rho$, resp. its restriction to the Cartan subalgebra corresponding to $\T\cap\M_{\qP}$ and $F_{\mu}$ denote an irreducible $\qL_{\qP}$-module of highest weight $\mu$. 
As a result, $$E_2^{p,q}=\bigoplus_{w\in W^{\qP}, \,l(w)=q}H^p(\mathfrak{l}_{\qP},K_{L_{\qP}};H_{\pi}\otimes F_{\mu_w}\otimes\mathrm{Sym}(\mathfrak{\check{a}}_{\qP}))$$
in view of which decomposition we give the 
\begin{defn}
A cohomology class in $H^*(\mathfrak{l}_{\qP},K_{L_{\qP}};H_{\pi}\otimes F_{\mu_w}\otimes\mathrm{Sym}(\mathfrak{\check{a}}_{\qP}))$ is called a \textit{class of type $(\pi,w)$}.
\end{defn}

\begin{rem}
 If the highest weight $\lambda_{\tau}$ of $(\tau,E)$ is regular, the weights $\mu_w$, $w\in W^{\qP}$ are regular  when restricted to the Cartan subalgebra corresponding to $\T\cap\M_{\qP}$.
\end{rem}

\section{Eisenstein cohomology in the generic case}\label{regular}

The highest weight $\lambda_{\tau}$ of the irreducible finite dimensional complex rational representation $(\tau,E)$ of $\G$ is now assumed to be regular. With notation as before the ad\`elic version of Theorem $4.11$ in \cite{Sch1} is

\begin{thm}(\cite{Sch1})\label{mainthm}
Let $\qP$ be a parabolic $\q$-subgroup of $\G$ and $\qQ\in\{\qP\}$ any element in its associate class.
If the Eisenstein series $E^{\G}_{\qQ}(f,\lambda)$ attached to a non-trivial cohomology class in $H^*(\mathfrak{g},K;W_{\qQ,\tilde{\pi}}\otimes\mathrm{Sym}(\mathfrak{\check{a}}_{\qQ})\otimes E)$ of type $(\pi,w)$, where $\pi\in\phi_{\qQ}$, $w\in W^{\qQ}$ and $f\in W_{\qQ,\tilde{\pi}}$,
is holomorphic at the point $\lambda_w=-w(\lambda_{\tau}+\rho)|_{\mathfrak{a}_{\qQ}}$, then $E^{\G}_{\qQ}(f,\lambda_w)$ represents a non-trivial cohomology class in $H^*(\mathfrak{g},K;\mathcal{A}_{E,\{\qP\},\phi}\otimes E)$.
\end{thm}

Concerning the question of holomorphy of the Eisenstein series $E^{\G}_{\qQ}(f,\lambda)$ at $\lambda_w$ an affirmative answer can be given in general for maximal parabolic $\q$-subgroups.

\begin{thm}(\cite{Sch4})\label{holomorphy}
Let $\mathbf{P}$ be a maximal parabolic $\q$-subgroup of $\G$ and $\qQ\in\{\qP\}$ any element in 
its associate class. Let $E^{\G}_{\qQ}(f,\lambda)$ be the Eisenstein series attached to a non-trivial cohomology class in $H^*(\mathfrak{g},K;W_{\qQ,\tilde{\pi}}\otimes\mathrm{Sym}(\mathfrak{\check{a}}_{\qQ})\otimes E)$ of 
type $(\pi,w)$, where $\pi\in\phi_{\qQ}$, $w\in W^{\qQ}$ such that $l(w)\geq\frac{1}{2}\mathrm{dim}_{\rea}\N_{\qQ}(\rea)$, and $f\in W_{\qQ,\tilde{\pi}}$. 
Then the meromorphic continuation of $E^{\G}_{\qQ}(f,\lambda)$ to $\mathfrak{\check{a}}_{\qQ}$ is holomorphic at $\lambda_w$.
\end{thm}

The vanishing theorem to follow was proven independently by J.-S. Li \& J. Schwermer and L. Saper and is cited from \cite{L-S}. 
The bounds for vanishing are given in terms of the two quantities: $l_0(\G(\rea))=\textrm{rk}(\G(\rea))-\textrm{rk}(K)$, where $\textrm{rk}$ is to denote the absolute rank of the group in question, and $q_0(\G(\rea))=\frac{1}{2}(\dim X-l_0(\G(\rea)))$.

\begin{thm}\label{vanishing}
If $\{\mathbf{P}\}\in\mathcal{C}$ is a class of associate parabolic $\q$-subgroups and $\{\mathbf{P}\}\neq\{\G\}$, then the summand $H^*(\mathfrak{g},K;\mathcal{A}_{E,\{\qP\}}\otimes E)$ 
is spanned by regular Eisenstein cohomology classes and 
$H^q(\mathfrak{g},K;\mathcal{A}_{E,\{\qP\}}\otimes E)=0$ for $q<q_0(\G(\rea))$.
\end{thm}

This is Theorem $5.5$ in \cite{L-S} for the case that the center of $\G$ is anisotropic over $\q$. 
\begin{rem}
Using Poincar\'e duality (cf. Theorem $5.6$ in \cite{L-S}) it can be inferred that in the generic case the cuspidal cohomology of $\G$ vanishes outside the interval $[q_0(\G),q_0(\G)+l_0(\G)]$.
\end{rem}

Finally, let $\textrm{vcd}(\G)=\dim X-\textrm{rk}_{\q}\G$ be the \textit{virtual cohomological dimension} of $\G$. Then the range of non-vanishing 
of the automorphic cohomology of $\G$ with respect to the degree is given by the interval $[q_0(\G(\rea)),\textrm{vcd}(\G)]$.

\section{A rational form of $SO(n,2)$}\label{q-form}

Let $V$ be a vector space over $\q$ of dimension $n+2$ and $f$ the regular quadratic form on $V$ 
of Witt-index $2$ represented by the symmetric matrix 
$$F=F_{n,2}=\left(\begin{array}{ccc}
0&0&I_2\\
0&I_{n-2}&0\\
I_2&0&0 
\end{array}\right)$$
with $I_n$ the $n\times n$ identity matrix. Let $\mathrm{SO}(f)$ denote the group of proper 
isometries of the rational quadratic space $(V,f)$ and $\G$ the $\q$-rational form of $\q$-rank $2$ of the self-adjoint semisimple linear algebraic group $\mathbf{SO}(n+2,\co)$ obtained from $\mathrm{SO}(f)$. 

The group $G=\G(\rea)$ of real points of $\G$ is isomorphic to the real semisimple Lie group $SO(n,2)$,
which has two connected components in the real topology. 

\begin{rem}
For the sake of a treatment of the groups $SO(n,2)$ -- with $n$ growing arbitrarily -- as uniform as possible, it is assumed from now on that $n\geq5$ implying, in particular, that $\G$ is not (quasi)-split.\footnote{For $n\in\{2,3,4\}$ the special orthogonal group $SO(n,2)$ is isogenous to (the direct product of) other semisimple real Lie groups, the real root system of which is of type 
$\textrm{A I}$, $\textrm{C I}$ and $\textrm{A III}$ in the Cartan numbering, respectively. As a matter of fact 
these cases are covered by the existing literature (cf. \cite{Har1}, \cite{Sch2}, \cite{Ha-Sch}). }
\end{rem}

\subsubsection*{The Lie algebra}

The Lie algebra $\mathfrak{g}$ of $\G(\co)$ is realised as the matrix algebra of 
complex $(n+2)\times(n+2)$-matrices $X$ being skew-symmetric with respect to $F$. 
Then the real form $\mathfrak{g}_0=\textrm{Lie}G$ is the subalgebra of all such $X$ with real entries. 
The Cartan involution on $\mathfrak{g}_0$ is written as $\theta:X\mapsto FXF$. 

As $\mathfrak{g}$ is isomorphic to $\mathfrak{so}(n+2,\co)$ the type of its root system depends on the parity of $n$. 
In particular, the type is $B_{\frac{n+1}{2}}$ if $2\nmid n$ and $D_{\frac{n+2}{2}}$ if $2\mid n$. 
Throughout, both cases will be treated separately, although an attempt is made to keep redundancy as little as possible.

\subsection{Rational roots and standard parabolic $\q$-subgroups $\qP$}
The subgroup $\tS=\{s_{a_1,a_2}=\textrm{diag}(a_1,a_2,1,...,1,a_1^{-1},a_2^{-1})\in \G\}$ is a 
maximal $\q$-split torus in $\G$ with rational characters $a_1:s_{a_1,a_2}\mapsto a_1$ and $a_2:s_{a_1,a_2}\mapsto a_2$. The relative root system of $(\G,\tS)$, which is of type $B_2$, is then given as the set $${}_{\q}\Phi=\{a_1a_2^{-1},a_1^{-1}a_2,a_1,a_2,a_1^{-1},a_2^{-1},a_1^{-1}a_2^{-1},a_1a_2\},$$ where 
moreover $\alpha_1=a_1a_2^{-1}$ and $\alpha_2=a_2$ are chosen to be simple. 

The parabolic $\q$-subgroups $\qP$ of $\G$ standard with respect to $\tS$ and ${}_{\q}\Phi^+$ are 
parametrised by the subsets $I$ of ${}_{\q}\Delta=\{\alpha_1,\alpha_2\}$. The maximal ones are listed below.\footnote{As for notation, the matrix $F_{k,l}$ with $k,l\in\mathbb{N}$ and possibly $l=0$ should always be interpreted in analogy to $F_{n,2}$ as at the beginning of the present section.}
In case of $I_1=\{\alpha_2\}$ the corresponding standard maximal parabolic $\q$-subgroup $\qP_1$ 
is described by
\small
$$
\begin{array}{l}
\A_1=\{s_{a_1,a_2}\in \tS|a_2=1\}\\
\\
\qL_1=\left\{g\in \G\Big\vert g=\left(\begin{array}{ccccc}
                           a_1&0&0&0&0\\
                           0&a_2&u_2&0&v_2\\
                           0&x_2&b&0&w_2\\
                           0&0&0&a_1^{-1}&0\\
                           0&y_2&z_2&0&c_2\end{array}\right)
\right\} \simeq \mathbf{SO}(F_{n-1,1})\times \mathbf{GL}_1\\
\\
\M_1=\left\{g\in \qL_1|a_1=\pm 1\right\} \simeq \mathbf{SO}(F_{n-1,1})\times\z/2\z\\
\\
\N_1=\left\{\left(\begin{array}{ccc}
                  a&u&v\\
                  0&I_{n-2}&-{}^tu\\
                  0&0&{}^ta^{-1} \end{array}\right)\Bigg\vert\begin{array}{l} a=\left(\begin{array}{cc}
                                              1&*\\
                                              0&1\end{array}\right), u=\left(\begin{array}{c}
                                              u_1\\
                                              0\end{array}\right),\\ 
v=a(r-\frac{1}{2}u{}^tu), r\in\mathbf{Mat}_{2\times2}, {}^tr=-r\end{array}\right\}\\
\\
\textrm{ is abelian and } \dim\N_1(\rea)=n.
\end{array}$$

\normalsize

In case of $I_2=\{\alpha_1\}$ the corresponding standard maximal parabolic $\q$-subgroup $\qP_2$ 
is described by
\small

$$
\begin{array}{l}
\A_2=\{s_{a_1,a_2}\in\tS|a_2=a_1\}\\
\\
\qL_2=\left\{g\in \G\Big\vert g=\left(\begin{array}{ccc}
                            a&0&0\\
                            0&b&0\\
                            0&0&{}^ta^{-1}\end{array}\right)
\right\} \simeq \mathbf{SO}(F_{n-2})\times \mathbf{GL}_2\\     
\\
\M_2=\left\{g\in \qL_2|a\in \mathbf{SL}_2^{\pm}\right\} \simeq \mathbf{SO}(F_{n-2})\times \mathbf{SL}_2^{\pm}\\
\\
\N_2=\left\{\left(\begin{array}{ccc}
                  I_2&u&r-\frac{1}{2}u{}^tu\\
                  0&I_{n-2}&-{}^tu\\
                  0&0&I_2
                 \end{array}\right)\Bigg\vert r\in\mathbf{Mat}_{2\times2}, {}^tr=-r\right\}\\
\\
\textrm{ is non-abelian and } \dim\N_2(\rea)=2n-3.
\end{array}$$
\normalsize

\begin{rem}
The conjugacy classes and the associate classes of parabolic $\q$-subgroups of $\G$ coincide.
\end{rem}

\subsection{The absolute root system}\label{absolute}
A maximal $\q$-torus $\T\subset\G$ containing $\tS$ is given by $\T=\tS\times\mathbf{SO}(F_2)^{k-2}$.
A positive system $\Phi^+$ compatible with ${}_{\q}\Phi^+$ is provided in terms of the Lie algebra $\mathfrak{h}$ of its centraliser $\qH=\T$.

\subsection*{The odd case}
Let $2\nmid n$ and set $2k=n+1$. The root system of $\mathfrak{g}$ is of type $B_k$ in this case.
In the realisation of $\mathfrak{g}$ mentioned above the Cartan subalgebra $\mathfrak{h}$ consists of complex block-diagonal matrices $H=\textrm{diag}(a,b,-a)$, where $a=\textrm{diag}(a_1,a_2)$ and $b$ itself is a skew-symmetric $2\times 2$-block diagonal matrix with the $(n-2)^{nd}$ column and row being zero. A block of $b$ is written in the form
$$\left(\begin{array}{cc}
0&b_j\\
-b_j&0\end{array}\right)$$ with $j=2i+1$ denoting the $j$-th row for $1\leq i\leq k-2$. 
The set of roots for $\mathfrak{h}$ in $\mathfrak{g}$ is 
$$\Phi=\Phi(\mathfrak{g},\mathfrak{h})=\{\pm\varepsilon_i\pm\varepsilon_j, \pm\varepsilon_i\vert1\leq i<j\leq k\}$$
with linear functionals $\varepsilon_1(H)=a_1$, $\varepsilon_2(H)=a_2$ on the symmetric $a$-block and $\varepsilon_{i+2}(H)=\sqrt{-1}b_j$ for $1\leq i\leq k-2$ on the skew-symmetric $b$-block.
The simple roots in $\Phi$ are chosen to be $\alpha_i=\varepsilon_i-\varepsilon_{i+1}$
for $1\leq i\leq k-1$ and $\alpha_k=\varepsilon_k$.
Thus, the restriction of the root $\alpha_1$, resp. $\alpha_2$ to $\mathfrak{s}$ is the differential at the identity of the character $\alpha_1$, resp. $\alpha_2$ in ${}_{\q}\Delta$. 

The standard maximal parabolic $\q$-subgroups $\qP_i$ are depicted via crossed Dynkin diagrams. 
$$\begin{array}{cccc}
\qP_1&&&\qP_2\\
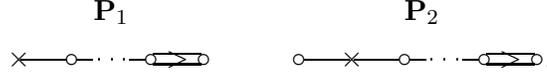
\begin{picture}(76,14)
\put(4,3){\line(1,0){17}}
\put(25,3){\line(1,0){6}}
\put(51,3){\line(-1,0){6}}
\put(54,1,2){\line(1,0){18}}
\put(54,5){\line(1,0){18}}
\put(39,3){\makebox(0,0){\dots}}
\put(63,3){\makebox(0,0){$>$}}
\put(3,3){\makebox(0,0){$\times$}}
\put(23,2,7){\makebox(0,0){$\circ$}}
\put(53,2,7){\makebox(0,0){$\circ$}}
\put(73,2,7){\makebox(0,0){$\circ$}}
\end{picture}&&&
\begin{picture}(96,14)
\put(5,3){\line(1,0){17}}
\put(24,3){\line(1,0){17}}
\put(45,3){\line(1,0){6}}
\put(71,3){\line(-1,0){6}}
\put(74,1,2){\line(1,0){18}}
\put(74,5){\line(1,0){18}}
\put(59,3){\makebox(0,0){\dots}}
\put(83,3){\makebox(0,0){$>$}}
\put(3,2,7){\makebox(0,0){$\circ$}}
\put(23,3){\makebox(0,0){$\times$}}
\put(43,2,7){\makebox(0,0){$\circ$}}
\put(73,2,7){\makebox(0,0){$\circ$}}
\put(93,2,7){\makebox(0,0){$\circ$}}
\end{picture}.\end{array}$$
The set $\Delta(\qL_i)$ of simple roots of the 
reductive algebraic $\q$-group $\qL_i$ can then be read off directly from the diagram as
$$ \Delta(\qL_1)=\Delta\setminus\{\alpha_1\},\quad \Delta(\qL_2)=\Delta\setminus\{\alpha_2\}.$$

\subsection*{The even case}
Let now $2\mid n$ and set $2k=n+2$. Then the root system of $\mathfrak{g}$ is of type $D_k$.
The Cartan subalgebra $\mathfrak{h}\subset\mathfrak{g}$ corresponding to the choice of $\T$ has the same form as in the odd case, this time, however, with $b$ being \textit{of maximal rank}.
With notation analogous to the odd case the set of roots for $\mathfrak{h}$ in $\mathfrak{g}$ is 
$$\Phi=\Phi(\mathfrak{g},\mathfrak{h})=\{\pm\varepsilon_i\pm\varepsilon_j\vert1\leq i<j\leq k\}$$ and   
the simple roots in $\Phi$ are chosen to be $\alpha_i=\varepsilon_i-\varepsilon_{i+1}$ for $1\leq i\leq k-1$ and $\alpha_k=\varepsilon_{k-1}+\varepsilon_k$.
The positive system $\Phi^+$ defined by the simple roots is compatible 
with the one on ${}_{\q}\Phi$ defined by ${}_{\q}\Delta$ and the crossed Dynkin diagrams for the standard maximal parabolic $\q$-subgroups $\qP_i$ are 
$$\begin{array}{cccc}
\qP_1&&&\qP_2\\
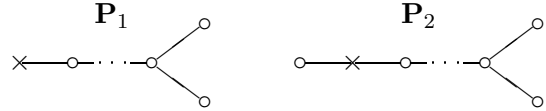
\begin{picture}(76,14)
\put(4,3){\line(1,0){17}}
\put(25,3){\line(1,0){6}}
\put(51,3){\line(-1,0){6}}
\put(55,4,5){\line(4,3){16}}
\put(55,1,5){\line(4,-3){16}}
\put(39,3){\makebox(0,0){\dots}}
\put(3,3){\makebox(0,0){$\times$}}
\put(23,2,7){\makebox(0,0){$\circ$}}
\put(53,2,7){\makebox(0,0){$\circ$}}
\put(73,-12){\makebox(0,0){$\circ$}}
\put(73,17,5){\makebox(0,0){$\circ$}}
\end{picture}&&&
\begin{picture}(96,14)
\put(5,3){\line(1,0){17}}
\put(24,3){\line(1,0){17}}
\put(45,3){\line(1,0){6}}
\put(71,3){\line(-1,0){6}}
\put(75,4,5){\line(4,3){16}}
\put(75,1,5){\line(4,-3){16}}
\put(59,3){\makebox(0,0){\dots}}
\put(3,2,7){\makebox(0,0){$\circ$}}
\put(23,3){\makebox(0,0){$\times$}}
\put(43,2,7){\makebox(0,0){$\circ$}}
\put(73,2,7){\makebox(0,0){$\circ$}}
\put(93,-12){\makebox(0,0){$\circ$}}
\put(93,17,5){\makebox(0,0){$\circ$}}
\end{picture}\end{array}$$

\vspace*{4mm}

\noindent which exhibit the sets $\Delta(\qL_i)$ to be formally the same as in the odd case.

\subsection{The Levi subalgebras}\label{Levis}

Consider the first standard maximal parabolic $\q$-subgroup $\qP_1$. Then $\mathfrak{m}_1=\textrm{Lie}(\M_1(\co))$ is realised relative to $F_{n-1,1}$. 
Set $\mathfrak{b}_1=\mathfrak{h}\cap\mathfrak{m}_1$ for the Cartan subalgebra of $\mathfrak{m}_1$ obtained 
from $\mathfrak{h}\subset\mathfrak{g}$, then
$$\mathfrak{b}_1=\left\{X\in\mathfrak{m}_1\vert X=\textrm{diag}(a_2,b,-a_2), a_2\in\co, \textrm{diag}(0,b,0)\in\mathfrak{h}\right\}.
$$

Consider now the second standard maximal parabolic $\q$-subgroup $\qP_2$. Then $\mathfrak{m}_2=\textrm{Lie}(\M_2(\co))$ is 
the Lie subalgebra of $\mathfrak{g}$ consisting of block diagonal matrices $\textrm{diag}(a,b,-a)\in\mathfrak{g}$
with $a$ of trace zero.
Set $\mathfrak{b}_2=\mathfrak{h}\cap\mathfrak{m}_2$ for the Cartan subalgebra of $\mathfrak{m}_2$ obtained 
from $\mathfrak{h}\subset\mathfrak{g}$, then
$$\mathfrak{b}_2=\left\{X\in\mathfrak{m}_2\vert a=\textrm{diag}(a_1,-a_1), \textrm{diag}(0,b,0)\in\mathfrak{h}\right\}.
$$

As for $i\in\{1,2\}$ the Levi subalgebra $\mathfrak{l}_i\cong\mathfrak{m}_i\oplus\mathfrak{a}_i$ with $\mathfrak{a}_i$ the center of $\mathfrak{l}_i$, the Cartan subalgebra $\mathfrak{h}$ decomposes as $\mathfrak{h}\cong\mathfrak{b}_i\oplus\mathfrak{a}_i$, which is even true in terms of the real forms.

\subsection{Fundamental weights and coordinates}\label{fundamental}
The fundamental 
weights defined by $\Delta$ are denoted in Bourbaki notation (cf. \cite{Bou}, Planche II, IV) and every weight $\lambda\in\mathfrak{\check{h}}$ will be written in coordinates $\lambda=(\lambda_1,...,\lambda_k)$ relative to the 
fundamental weights. Also, a dominant weight $\lambda$ is regular if and only if $\lambda_i>0$ for all $i\in\{1,\ldots,k\}$. 
Suppose $\mathfrak{\check{a}}_i$, resp. $\mathfrak{\check{b}}_i$ be identified with the space of all linear forms 
on $\mathfrak{h}$ vanishing on $\mathfrak{b}_{i}$, resp. $\mathfrak{a}_{i}$. In this way a canonical 
isomorphism $\mathfrak{\check{h}}\cong\mathfrak{\check{a}}_i\oplus\mathfrak{\check{b}}_i$ can be obtained, which in turn 
allows to restrict weights on $\mathfrak{h}$ in a canonical way to its direct summands. 

\subsection*{The odd case}
Let $2\nmid n$ and set $2k=n+1$. A basis of fundamental weights for each of the direct summands of the Cartan subalgebra is given in case of $i=1$ by
\small
$$\begin{array}{l}
\varpi_{11}=\varpi_1\in\mathfrak{\check{a}}_{1} \quad\mathrm{and}\\
\varpi_{12}=(-1,1,0,\dots,0), \varpi_{13}=(-1,0,1,0,\dots,0),\dots\\
\ldots,\varpi_{1k}=(-\frac{1}{2},0,\dots,0,1)
\in\mathfrak{\check{b}}_{1}\end{array}$$
\normalsize
and in case $i=2$ by
\small
$$\begin{array}{l}
\varpi_{21}=\varpi_2\in\mathfrak{\check{a}}_{2} \quad\mathrm{and}\\
\varpi_{22}=(1,-\frac{1}{2},0,\dots,0), \varpi_{23}=(0,-1,1,0,\dots,0),\dots\\
\ldots, \varpi_{2k}=(0,-\frac{1}{2},0,\dots,0,1)
\in\mathfrak{\check{b}}_{2}.\end{array}$$
\normalsize
Accordingly, the restriction of $\lambda=(\lambda_1,\dots,\lambda_k)\in\mathfrak{\check{h}}$ to $\mathfrak{a}_{i}$ and $\mathfrak{b}_{i}$ is given by:
$$\begin{array}{l}
(\lambda_1,\dots,\lambda_{k})|_{\mathfrak{a}_{1}}=(\lambda_1+\dots+\lambda_{k-1}+\frac{1}{2}\lambda_k)\varpi_{11}\\
(\lambda_1,\dots,\lambda_k)|_{\mathfrak{b}_{1}}=\lambda_2\varpi_{12}+\dots+\lambda_k\varpi_{1k}\end{array}$$
and
$$\begin{array}{l}
(\lambda_1,\dots,\lambda_k)|_{\mathfrak{a}_{2}}=
(\frac{1}{2}\lambda_1+\lambda_2+\dots+\lambda_{k-1}+\frac{1}{2}\lambda_k)\varpi_{21}\\
(\lambda_1,\dots,\lambda_k)|_{\mathfrak{b}_{2}}=\lambda_1\varpi_{22}+\lambda_3\varpi_{23}+\dots+\lambda_k\varpi_{2k}
\end{array}$$

Finally, let $\rho$ denote half the sum of the positive roots in $\Phi$ and $\rho_i=\rho|_{\mathfrak{a}_i}$ its restriction to $\mathfrak{a}_i$, $i=1,2$.
Then their coordinates are
$$\begin{array}{c}
\rho=(1,\dots,1),\quad \rho_1=(\frac{n}{2},0,\dots,0), \quad\rho_2=(0,\frac{n-1}{2},0,\dots,0).\end{array}$$

\subsection*{The even case}

Let $2\mid n$ and set $2k=n+2$. A basis of fundamental weights for each of the direct summands of the Cartan subalgebra in case of $i=1$:
\small
$$\begin{array}{l}
\varpi_{11}=\varpi_1\in\mathfrak{\check{a}}_{1} \textrm{ and}\\
\varpi_{12}=(-1,1,0,\dots,0),\dots, \varpi_{1k-2}=(-1,0\ldots,0,1,0,0),\\
\varpi_{1k-1}=(-\frac{1}{2},0,\dots,0,1,0),\varpi_{1k}=(-\frac{1}{2},0,\dots,0,1)
\in\mathfrak{\check{b}}_{1}\end{array}$$
\normalsize
and in case $i=2$:
\small
$$\begin{array}{l}
\varpi_{21}=\varpi_2\in\mathfrak{\check{a}}_{2} \textrm{ and}\\
\varpi_{22}=(1,-\frac{1}{2},0,\dots,0), \varpi_{23}=(0,-1,1,0,\dots,0),\dots\\
\ldots, \varpi_{2k-2}=(0,-1,0\ldots,0,1,0,0),\varpi_{2k-1}=(0,-\frac{1}{2},0,\dots,0,1,0),\\
\varpi_{2k}=(0,-\frac{1}{2},0,\dots,0,1)
\in\mathfrak{\check{b}}_{2}.\end{array}$$
\normalsize
Accordingly, the restriction of $\lambda=(\lambda_1,\dots,\lambda_k)\in\mathfrak{\check{h}}$ to 
$\mathfrak{a}_{i}$ and $\mathfrak{b}_{i}$ is given by:
$$\begin{array}{l}
(\lambda_1,\dots,\lambda_{k})|_{\mathfrak{a}_{1}}=(\lambda_1+\dots+\frac{1}{2}\lambda_{k-1}+\frac{1}{2}\lambda_k)
\varpi_{11}\\
(\lambda_1,\dots,\lambda_k)|_{\mathfrak{b}_{1}}=\lambda_2\varpi_{12}+\dots+\lambda_k\varpi_{1k}\end{array}$$
and
$$\begin{array}{l}
(\lambda_1,\dots,\lambda_k)|_{\mathfrak{a}_{2}}=
(\frac{1}{2}\lambda_1+\lambda_2+\dots+\frac{1}{2}\lambda_{k-1}+\frac{1}{2}\lambda_k)\varpi_{21}\\
(\lambda_1,\dots,\lambda_k)|_{\mathfrak{b}_{2}}=\lambda_1\varpi_{22}+\lambda_3\varpi_{23}+\dots+\lambda_k\varpi_{2k}
\end{array}$$
The coordinates of $\rho$ and $\rho_i=\rho|_{\mathfrak{a}_i}$ are formally identical to the ones in the odd case.

\section{Classes of type $(\pi,w)$}\label{oftype}

The cohomological contribution of the summand $H^q(\mathfrak{g},K;\mathcal{A}_{E,\{\qP\}}\otimes E)$,
$\qP$ a proper standard parabolic $\q$-subgroup of $\G$, to the automorphic cohomology of $\G$ can be reconstructed from cuspidal classes of types $(\pi,w)$. Their determination starts with the second parameter. As for the question of how to obtain the elements of $W^{\qP}$ explicitly 
appendix \ref{diagram} establishes an algorithm based on the results of \cite{Kos}.

\subsection{Minimal coset representatives}\label{Kostant}

The Weyl group $W=W(\Phi)$ of the root system $\Phi=\Phi(\G,\T)=\Phi(\mathfrak{g},\mathfrak{h})$
is generated by the reflections $s_{\alpha_i}$ at the simple roots $\alpha_i$ for $1\leq i\leq k$
from section \ref{absolute}. For the sake of convenience the generators of $W$ will be written as 
$\{s_1,\ldots,s_k\}$. 

The following sets out to display the elements of $W^{\qP}$ for the standard maximal parabolic
$\q$-subgroups of $\G$ uniformly for all $n\geq5$. Along with $w\in W^{\qP}$ the number $N(l)$ of such per length $l=l(w)$ will be provided. 

\subsubsection{The first standard maximal parabolic $\q$-subgroup}
The Weyl group $W_{\qL_1}$ is generated by the elements $s_2,\ldots,s_k$ given that the positive system on $\Phi(\qL_1,\T)$ is the one induced by $\Phi^+$. Let $W^{\qP_1}$ denote the set of minimal coset representatives of $W_{\qL_1}\backslash W$.

\subsection*{The odd case} 
Let $2\nmid n$ and set again $2k=n+1$. The order of $W^{\qP_1}$ is given as $|W^{\qP_1}|=n+1$ and the set $W^{\qP_1}$ is actually computed in table $1$.

\small
\begin{table}[ht]
\caption{$w\in W^{\qP_1}$}
\begin{tabular}{c l c}
$l(w)$ & $w$ & $N(l)$\\
\hline

$0$ & $1$ & $1$\\
$1$ & $s_1$ & $1$\\
$2$ & $s_1s_2$ & $1$\\
\vdots & \;\vdots & \vdots\\
$k$ & $s_1\cdots s_{k-1}s_k$ & $1$\\
$k+1$ & $s_1\cdots s_{k-1}s_ks_{k-1}$ & $1$\\
\vdots & \;\vdots & \vdots\\
$n-2$ & $s_1\cdots s_{k-1}s_ks_{k-1}\cdots s_3$ & $1$\\
$n-1$ & $s_1\cdots s_{k-1}s_ks_{k-1}\cdots s_3s_2$ & $1$\\
$n$ & $s_1\cdots s_{k-1}s_ks_{k-1}\cdots s_3s_2s_1$ & $1$\\
\end{tabular}

\end{table}
\normalsize

\subsection*{The even case}

Let $2\mid n$ and $2k=n+2$. 
The order of $W^{\qP_1}$ is given as $|W^{\qP_1}|=n+2$ and the set $W^{\qP_1}$ is provided by 
table $2$.

\small
\begin{table}[ht]
\caption{$w\in W^{\qP_1}$}
\begin{tabular}{c l c}
$l(w)$ & $w$ & $N(l)$\\
\hline

$0$ & $1$ & $1$\\
$1$ & $s_1$ & $1$\\
$2$ & $s_1s_2$ & $1$\\
\vdots & \;\vdots & \vdots\\
$k-2$ & $s_1\cdots s_{k-2}$ & $1$\\
$k-1$ & $\begin{cases}
          s_1\cdots s_{k-2}s_{k-1}\\s_1\cdots s_{k-2}s_{k}
         \end{cases}$
& $2$\\
$k$ & $s_1\cdots s_{k-2}s_{k-1}s_{k}$ & $1$\\
\vdots & \;\vdots & \vdots\\
$n-2$ & $s_1\cdots s_{k-2}s_{k-1}s_{k}s_{k-2}\cdots s_3$ & $1$\\
$n-1$ & $s_1\cdots s_{k-2}s_{k-1}s_{k}s_{k-2}\cdots s_3s_2$ & $1$\\
$n$ & $s_1\cdots s_{k-2}s_{k-1}s_{k}s_{k-2}\cdots s_3s_2s_1$ & $1$\\
\end{tabular}

\end{table}
\normalsize

\subsubsection{The second standard maximal parabolic $\q$-subgroup}
The Weyl group $W_{\qL_2}$ is generated by the elements $s_1,s_3,\ldots,s_k$.
Let $W^{\qP_2}$ denote the set of minimal coset representatives of $W_{\qL_2}\backslash W$. 

\subsection*{The odd case} 
Let $2\nmid n$ and set $2k=n+1$. The order of $W^{\qP_2}$ is given as $|W^{\qP_2}|=\frac{(n+1)(n-1)}{2}$ and the set $W^{\qP_2}$ is computed in table $3$.

\small
\begin{table}[ht]
\caption{$w\in W^{\qP_2}$}
\begin{tabular}{c l c}
$l(w)$ & $w$ & $N(l)$\\
\hline

$0$ & $1$ & $1$\\
$1$ & $s_2$ & $1$\\
$2$ & $\begin{cases}
        s_2s_3\\s_2s_1
       \end{cases}$ & $2$\\
$3$ & $\begin{cases}
        s_2s_3s_4\\s_2s_1s_3
       \end{cases}$ & $2$\\
\vdots & \;\vdots & \vdots\\
$n-3$ & $\begin{cases}
          s_2s_3\cdots s_k\cdots s_3\mathbf{,} s_2s_1s_3s_4\cdots s_k\cdots s_4\mathbf{,}\\
	  s_2s_1s_3s_2s_4\cdots s_k\cdots s_5\mathbf{,}\ldots\\
	  \ldots\mathbf{,}\; s_2s_1s_3s_2s_4s_3\cdots s_{k-1}s_{k-2}
         \end{cases}$ & $\frac{n-1}{2}$\\
$n-2$ & $\begin{cases}
          s_2s_3\cdots s_k\cdots s_3s_2\mathbf{,} s_2s_1s_3s_4\cdots s_k\cdots s_4s_3\mathbf{,}\\
	  s_2s_1s_3s_2s_4\cdots s_k\cdots s_5s_4\mathbf{,}\ldots\\
	  \ldots\mathbf{,}\; s_2s_1s_3s_2s_4s_3\cdots s_{k-1}s_{k-2}s_k
         \end{cases}$ & $\frac{n-1}{2}$\\
$n-1$ & $\begin{cases}
          s_2s_3\cdots s_k\cdots s_2s_1\mathbf{,} s_2s_1s_3s_4\cdots s_k\cdots s_3s_2\mathbf{,}\\
	  s_2s_1s_3s_2s_4\cdots s_k\cdots s_4s_3\mathbf{,}\ldots\\
	  \ldots\mathbf{,}\; s_2s_1s_3s_2s_4s_3\cdots 
	  s_{k-1}s_{k-2}s_ks_{k-1}
         \end{cases}$ & $\frac{n-1}{2}$\\
$n$ & $\begin{cases}
          s_2s_3\cdots s_k\cdots s_2s_1s_2\mathbf{,} s_2s_1s_3s_4\cdots s_k\cdots s_3s_2s_3\mathbf{,}\\
	  s_2s_1s_3s_2s_4\cdots s_k\cdots s_4s_3s_4\mathbf{,}\ldots\\
	  \ldots\mathbf{,}\; s_2s_1s_3s_2s_4s_3\cdots 
	  s_{k-1}s_{k-2}s_ks_{k-1}s_k
         \end{cases}$ & $\frac{n-1}{2}$\\

\vdots & \;\vdots & \vdots\\
$2n-6$ & $\begin{cases}
           s_2s_3\cdots s_k\cdots s_2s_1s_2 \cdots s_k\cdots s_5\\
	   s_2s_1s_3s_4\cdots s_k\cdots s_3s_2s_3 \cdots s_k\cdots s_4
          \end{cases}$ & $2$\\
$2n-5$ & $\begin{cases}
           s_2s_3\cdots s_k\cdots s_2s_1s_2 \cdots s_k\cdots s_5s_4\\
	   s_2s_1s_3s_4\cdots s_k\cdots s_3s_2s_3 \cdots s_k\cdots s_4s_3
          \end{cases}$ & $2$\\
$2n-4$ & $s_2s_3\cdots s_k\cdots s_2s_1s_2 \cdots s_k\cdots s_5s_4s_3$ & $1$ \\
$2n-3$ & $s_2s_3\cdots s_k\cdots s_2s_1s_2 \cdots s_k\cdots s_5s_4s_3s_2$ & $1$ 

\end{tabular}
\end{table}

\normalsize

\begin{rem}
As it stands table $3$ applies to cases $n\geq9$. For $n=5,7$ obvious minor 
adjustments are in place. 
\end{rem}

\subsection*{The even case}

Let $2\mid n$ and $2k=n+2$. The set $W^{\qP_2}$ has order $|W^{\qP_2}|=\frac{(n+2)n}{2}$ and its elements are displayed in table $4$.

\small
\begin{table}[ht]
\caption{$w\in W^{\qP_2}$}
\begin{tabular}{c l c}
$l(w)$ & $w$ & $N(l)$\\
\hline

$0$ & $1$ & $1$\\
$1$ & $s_2$ & $1$\\
$2$ & $\begin{cases}
        s_2s_3\\s_2s_1
       \end{cases}$ & $2$\\
$3$ & $\begin{cases}
        s_2s_3s_4\\s_2s_1s_3
       \end{cases}$ & $2$\\
\vdots & \;\vdots & \vdots\\
$n-3$ & $\begin{cases}
         s_2s_3\cdots s_ks_{k-2}\cdots s_3\mathbf{,} s_2s_1s_3s_4\cdots s_ks_{k-2}\cdots s_4\\
	 s_2s_1s_3s_2s_4s_5\cdots s_ks_{k-2}\cdots s_5 \mathbf{,}\ldots\\
	 \ldots\mathbf{,}\; s_2s_1s_3s_2s_4s_3\cdots s_{k-3}s_{k-4}s_{k-2}s_{k-1}s_k\mathbf{,}\\
	 s_2s_1s_3s_2s_4s_3\cdots s_{k-2}s_{k-3}\begin{cases}s_{k-1}\\s_k\end{cases}
        \end{cases}$ & $\frac{n}{2}$\\
$n-2$ & $\begin{cases}
         s_2s_3\cdots s_ks_{k-2}\cdots s_3s_2\mathbf{,} s_2s_1s_3s_4\cdots s_ks_{k-2}\cdots s_4s_3\\
	 s_2s_1s_3s_2s_4s_5\cdots s_ks_{k-2}\cdots s_5s_4 \mathbf{,}\ldots\\
	 \ldots\mathbf{,}\; s_2s_1s_3s_2s_4s_3\cdots s_{k-3}s_{k-4}s_{k-2}s_{k-1}s_ks_{k-2}\mathbf{,}\\
	 s_2s_1s_3s_2s_4s_3\cdots s_{k-2}s_{k-3}\begin{cases}s_{k-1}s_k\\s_ks_{k-2}\\s_{k-1}s_{k-2}\end{cases}
        \end{cases}$ & $\frac{n+2}{2}$\\
$n-1$ & $\begin{cases}
         s_2s_3\cdots s_ks_{k-2}\cdots s_2s_1\mathbf{,} s_2s_1s_3s_4\cdots s_ks_{k-2}\cdots s_3s_2\\
	 s_2s_1s_3s_2s_4\cdots s_ks_{k-2}\cdots s_4s_3 \mathbf{,}\ldots\\
	 \ldots\mathbf{,}\; s_2s_1s_3s_2s_4s_3\cdots s_{k-3}s_{k-4}s_{k-2}s_{k-1}s_ks_{k-2}s_{k-3}\mathbf{,}\\
	 s_2s_1s_3s_2s_4s_3\cdots 
	 s_{k-2}s_{k-3}\begin{cases}s_{k-1}s_ks_{k-2}\\s_ks_{k-2}s_{k-1}\\s_{k-1}s_{k-2}s_k\end{cases}
        \end{cases}$ & $\frac{n+2}{2}$\\
$n$ & $\begin{cases}
         s_2s_3\cdots s_ks_{k-2}\cdots s_2s_1s_2\mathbf{,} s_2s_1s_3s_4\cdots s_ks_{k-2}\cdots s_3s_2s_3\\
	 s_2s_1s_3s_2s_4\cdots s_ks_{k-2}\cdots s_4s_3s_4 \mathbf{,}\ldots\\
	 \ldots\mathbf{,}\; s_2s_1s_3s_2s_4s_3\cdots s_{k-3}s_{k-4}s_{k-2}s_{k-1}s_ks_{k-2}s_{k-3}s_{k-2}\mathbf{,}\\
	 s_2s_1s_3s_2s_4s_3\cdots 
	 s_{k-2}s_{k-3}s_{k-1}s_ks_{k-2}\begin{cases}s_k\\s_{k-1}\end{cases}
        \end{cases}$ & $\frac{n}{2}$\\
\vdots & \;\vdots & \vdots\\
$2n-6$ & $\begin{cases}
           s_2s_3\cdots s_ks_{k-2}\cdots s_2s_1s_2 \cdots s_{k-2}s_k\cdots s_5\\
	   s_2s_1s_3s_4\cdots s_ks_{k-2}\cdots s_3s_2s_3 \cdots s_{k-2}s_k\cdots s_4
          \end{cases}$ & $2$\\
$2n-5$ & $\begin{cases}
           s_2s_3\cdots s_ks_{k-2}\cdots s_2s_1s_2 \cdots s_{k-2}s_k\cdots s_5s_4\\
	   s_2s_1s_3s_4\cdots s_ks_{k-2}\cdots s_3s_2s_3 \cdots s_{k-2}s_k\cdots s_4s_3
          \end{cases}$ & $2$\\
$2n-4$ & $s_2s_3\cdots s_ks_{k-2}\cdots s_2s_1s_2 \cdots s_{k-2}s_k\cdots s_5s_4s_3$ & $1$ \\
$2n-3$ & $s_2s_3\cdots s_ks_{k-2}\cdots s_2s_1s_2 \cdots s_{k-2}s_k\cdots s_5s_4s_3s_2$ & $1$ 

\end{tabular}
\end{table}

\normalsize

\begin{rem}
As it stands table $4$ is valid for cases $n\geq10$ only. Again the 
remaining cases $n=6,8$ can be easily supplied by minor adjustments. 
\end{rem}

\subsection{The simple $\qL_{\qP}$-modules of Kostant's theorem}\label{Lmodules}

As an immediate application of the determination of the minimal coset representatives this section provides 
a listing of weights $\mu_w=w(\lambda+\rho)-\rho$ for $w\in W^{\qP_i}$, resp. their restrictions to the 
Cartan subalgebra $\mathfrak{b}_i$ for $i=1,2$ occurring in Kostant's theorem. Notation is as in 
\ref{fundamental}.

\subsubsection*{The first standard maximal parabolic $\q$-subgroup}

Let $\lambda\in\mathfrak{\check{h}}$ be a weight. The restriction of $\lambda=(\lambda_1,\ldots,\lambda_k)$ to the Cartan subalgebra $\mathfrak{b}_1$ is formally identical in the odd and even case: 
$(\lambda_1,\dots,\lambda_k)|_{\mathfrak{b}_{1}}=\lambda_2\varpi_{12}+\dots+\lambda_k\varpi_{1k}$.

\subsection*{The odd case}
Let $2\nmid n$ and set $2k=n+1$. Table $5$ lists the weights $w(\lambda+\rho)-\rho_{|\mathfrak{b}_1}$
in this case.

\scriptsize
\begin{table}[ht]
\caption{Restriction of $\mu_w=w(\lambda+\rho)-\rho$ to $\mathfrak{b}_1$ for $w\in W^{\qP_1}$: odd case}
\begin{tabular}{c l l}
$l(w)$ & $w$ & ${\mu_w}_{|\mathfrak{b}_1}$\\
\hline

$0$ & $1$ & $(\lambda_2,\ldots,\lambda_k)$\\
$1$ & $s_1$ & $(\lambda_2+\lambda_1+1,\lambda_3,\ldots,\lambda_k)$\\
\vdots & \;\vdots & \quad\vdots\\
$i$ & $s_1\cdots s_i$ & $(\lambda_1,\ldots,\lambda_{i-1},\lambda_{i+1}+\lambda_i+1,\lambda_{i+2},\ldots,\lambda_k)$\\
\vdots & \;\vdots & \quad\vdots\\
$k-1$ & $s_1\cdots s_{k-1}$ & $(\lambda_1,\ldots,\lambda_{k-2},\lambda_k+2\lambda_{k-1}+2)$\\
$k$ & $s_1\cdots s_{k-1}s_k$ & $(\lambda_1,\ldots,\lambda_{k-2},\lambda_k+2\lambda_{k-1}+2)$\\
\vdots & \;\vdots & \quad\vdots\\
$n-i$ & $s_1\cdots s_{k-1}s_ks_{k-1}\cdots s_{i+1}$ & $(\lambda_1,\ldots,\lambda_{i-1},\lambda_{i+1}+\lambda_i+1,\lambda_{i+2},\ldots,\lambda_k)$\\
\vdots & \;\vdots & \quad\vdots\\
$n-1$ & $s_1\cdots s_{k-1}s_ks_{k-1}\cdots s_2$ & $(\lambda_2+\lambda_1+1,\lambda_3,\ldots,\lambda_k)$\\
$n$ & $s_1\cdots s_{k-1}s_ks_{k-1}\cdots s_2s_1$ & $(\lambda_2,\ldots,\lambda_k)$
\end{tabular}

\end{table}
\normalsize

\begin{rem}
For reasons of consistency it is assumed that $2\leq i\leq k-2$ in table $5$.
\end{rem}

\subsection*{The even case}
Let $2\mid n$ and set $2k=n+2$. For table $6$ below it is assumed that $2\leq i\leq k-3$.

\scriptsize
\begin{table}[ht]
\caption{Restriction of $\mu_w=w(\lambda+\rho)-\rho$ to $\mathfrak{b}_1$ for $w\in W^{\qP_1}$: even case}
\begin{tabular}{c l l}
$l(w)$ & $w$ & ${\mu_w}_{|\mathfrak{b}_1}$\\
\hline

$0$ & $1$ & $(\lambda_2,\ldots,\lambda_k)$\\
$1$ & $s_1$ & $(\lambda_2+\lambda_1+1,\lambda_3,\ldots,\lambda_k)$\\
\vdots & \;\vdots & \quad\vdots\\
$i$ & $s_1\cdots s_i$ & $(\lambda_1,\ldots,\lambda_{i-1},\lambda_{i+1}+\lambda_i+1,\lambda_{i+2},\ldots,\lambda_k)$\\
\vdots & \;\vdots & \quad\vdots\\
$k-2$ & $s_1\cdots s_{k-2}$ & $(\lambda_1,\ldots,\lambda_{k-3},\lambda_{k-1}+\lambda_{k-2}+1,\lambda_k+\lambda_{k-2}+1)$\\
$k-1$ & $\begin{cases}s_1\cdots s_{k-2}s_{k-1}\\s_1\cdots s_{k-2}s_k          
         \end{cases}$ & 
	 $\begin{array}{c}
 (\lambda_1,\ldots,\lambda_{k-3},\lambda_{k-2},\lambda_k+\lambda_{k-2}+\lambda_{k-1}+2)\\
 (\lambda_1,\ldots,\lambda_{k-3},\lambda_k+\lambda_{k-2}+\lambda_{k-1}+2,\lambda_{k-2})
           \end{array}$\\
$k$ & $s_1\cdots s_{k-2}s_{k-1}s_k$ & $(\lambda_1,\ldots,\lambda_{k-3},\lambda_k+\lambda_{k-2}+1,\lambda_{k-1}+\lambda_{k-2}+1)$\\
\vdots & \;\vdots & \quad\vdots\\
$n-i$ & $s_1\cdots s_{k-1}s_ks_{k-1}\cdots s_{i+1}$ & $(\lambda_1,\ldots,\lambda_{i-1},\lambda_{i+1}+\lambda_i+1,\lambda_{i+2},\ldots,\lambda_k,\lambda_{k-1})$\\
\vdots & \;\vdots & \quad\vdots\\
$n-1$ & $s_1\cdots s_{k-1}s_ks_{k-1}\cdots s_2$ & $(\lambda_2+\lambda_1+1,\lambda_3,\ldots,\lambda_k,\lambda_{k-1})$\\
$n$ & $s_1\cdots s_{k-1}s_ks_{k-1}\cdots s_2s_1$ & $(\lambda_2,\ldots,\lambda_k,\lambda_{k-1})$
\end{tabular}

\end{table}
\normalsize

\subsubsection*{The second standard maximal parabolic $\q$-subgroup}
The restriction of $\lambda=(\lambda_1,\dots,\lambda_k)\in\mathfrak{\check{h}}$ to 
$\mathfrak{b}_{2}$ is formally identical in the odd and even case. It is given as
$(\lambda_1,\dots,\lambda_k)|_{\mathfrak{b}_{2}}=\lambda_1\varpi_{22}+\lambda_3\varpi_{23}+\dots+\lambda_k \varpi_{2k}$.

As always, the natural number $k$ is chosen such that $2k=n+1$, if $2\nmid n$, resp. $2k=n+2$, if $2\mid n$. Tables $7$ and $8$ provide the weights ${\mu_w}_{|\mathfrak{b}_2}$ in respective cases.

\begin{landscape}

\scriptsize
\begin{table}[ht]
\caption{Restriction of $\mu_w=w(\lambda+\rho)-\rho$ to $\mathfrak{b}_2$ for $w\in W^{\qP_2}$: odd case}
\begin{tabular}{c l l}
$l(w)$ & $w$ & ${\mu_w}_{|\mathfrak{b}_2}$\\
\hline

$0$ & $1$ & $(\lambda_1,\lambda_3,\ldots,\lambda_{k})$\\
$1$ & $s_2$ & $(\lambda_1+\lambda_2+1,\lambda_3+\lambda_2+1,\lambda_4,\ldots,\lambda_k)$\\
\vdots & \;\vdots & \quad\vdots\\
$k-1$ & $\begin{cases}s_2\cdots s_{k}\\s_2s_1s_3s_4\cdots s_{k-1}\\
\textrm{etc.} 
          \end{cases}$
      & $\begin{array}{l}
          (\lambda_1+\ldots+\lambda_{k}+(k-1),\lambda_2,\ldots,\lambda_{k-2},\lambda_k+2\lambda_{k-1}+2)\\
	  (\lambda_2+\ldots+\lambda_{k-1}+(k-3),\lambda_2+\lambda_1+1,\lambda_3,\ldots,\lambda_{k-2},
	  \lambda_k+2\lambda_{k-1}+2)\\
	  \textrm{etc.}
         \end{array}$\\
\vdots & \;\vdots & \quad\vdots \\
$n-2$ & $\begin{cases}s_2\cdots s_k\cdots s_2\\s_2s_1s_3s_4\cdots s_k\cdots s_4s_3\\
\vdots\\s_2s_1s_3s_2s_4s_3\cdots s_{k-1}s_{k-2}s_k          
         \end{cases}$ &
          $\begin{array}{l}
          (\lambda_1+\ldots+\lambda_k+\ldots+\lambda_2+(2k-3),\lambda_3,\ldots,\lambda_k)\\
	  (\lambda_2+\ldots+\lambda_k+\ldots+\lambda_3+(2k-5),\lambda_3+\lambda_2+\lambda_1+2,\lambda_4,
	  \ldots,\lambda_k)\\
                         \quad\vdots\\(\lambda_{k-1}+\lambda_k+1,\lambda_1,\ldots,\lambda_{k-3},\lambda_k+2(\lambda_{k-1}+\lambda_{k-2})+4)
	           \end{array}$\\
$n-1$ & $\begin{cases}s_2\cdots s_k\cdots s_2s_1\\s_2s_1s_3s_4\cdots s_k\cdots s_4s_3s_2\\
\vdots\\s_2s_1s_3s_2s_4s_3\cdots s_{k-2}s_ks_{k-1}       
         \end{cases}$ &
          $\begin{array}{l}
          (\lambda_1+\ldots+\lambda_k+\ldots+\lambda_2+(2k-3),\lambda_3,\ldots,\lambda_k)\\
	  (\lambda_2+\ldots+\lambda_k+\ldots+\lambda_3+(2k-5),\lambda_3+\lambda_2+\lambda_1+2,\lambda_4,
	  \ldots,\lambda_k)\\
\quad\vdots\\(\lambda_{k-1}+\lambda_k+1,\lambda_1,\ldots,\lambda_{k-3},\lambda_k+2(
	  \lambda_{k-1}+\lambda_{k-2})+4)
	           \end{array}$\\
\vdots & \;\vdots & \quad\vdots\\
$3k-4$ & $\begin{cases}s_2\cdots s_k\cdots s_2s_1s_2\cdots s_{k-1}\\s_2s_1s_3s_4\cdots s_k\cdots s_2\cdots s_k\\
\textrm{etc.} 
          \end{cases}$
      & $\begin{array}{l}
          (\lambda_1+\ldots+\lambda_{k}+(k-1),\lambda_2,\ldots,\lambda_{k-2},\lambda_k+2\lambda_{k-1}+2)\\
	 (\lambda_2+\ldots+\lambda_{k-1}+(k-3),\lambda_2+\lambda_1+1,\lambda_3,\ldots,\lambda_{k-2},
	  \lambda_k+2\lambda_{k-1}+2)\\
	  \textrm{etc.}
         \end{array}$\\
\vdots & \;\vdots & \quad\vdots\\
$2n-4$ & $s_2\dots s_k\dots s_2s_1s_2\dots s_k\dots s_3$ & 
$(\lambda_1+\lambda_2+1,\lambda_3+\lambda_2+1,\lambda_4,\ldots,\lambda_k)$\\
$2n-3$ & $s_2\cdots s_k\cdots s_2s_1s_2\cdots s_k\cdots s_2$ & $(\lambda_1,\lambda_3,\ldots,\lambda_{k})$
\end{tabular}

\end{table}
\normalsize

\end{landscape}

\begin{landscape}

\scriptsize
\begin{table}[ht]
\caption{Restriction of $\mu_w=w(\lambda+\rho)-\rho$ to $\mathfrak{b}_2$ for $w\in W^{\qP_2}$: even case}
\begin{tabular}{c l l}
$l(w)$ & $w$ & ${\mu_w}_{|\mathfrak{b}_2}$\\
\hline

$0$ & $1$ & $(\lambda_1,\lambda_3,\ldots,\lambda_{k})$\\
$1$ & $s_2$ & $(\lambda_1+\lambda_2+1,\lambda_3+\lambda_2+1,\lambda_4,\ldots,\lambda_k)$\\
\vdots & \;\vdots & \quad\vdots\\
$k-1$ & $\begin{cases}s_2\cdots s_k\\s_2s_1s_3s_4\cdots s_{k-1}\\
\textrm{etc.}\end{cases}$ & 
$\begin{array}{l}(\lambda_1+\ldots+\lambda_k+(k-1),\lambda_2,\ldots,\lambda_{k-3},\lambda_k+\lambda_{k-2}+1,\lambda_{k-1}+\lambda_{k-2}+1)\\(\lambda_2+\ldots+\lambda_{k-1}+(k-3),\lambda_2+\lambda_1+1,\lambda_3,\ldots,\lambda_{k-2},\lambda_k+\lambda_{k-1}+\lambda_{k-2}+2)\\
\textrm{etc.}\end{array}$\\
\vdots & \;\vdots & \quad\vdots\\
$n-2$ & $\begin{cases}
         s_2s_3\cdots s_ks_{k-2}\cdots s_3s_2\\s_2s_1s_3s_4\cdots s_ks_{k-2}\cdots s_4s_3\\
\vdots\\
	 s_2s_1s_3s_2s_4\cdots s_{k-2}s_{k-3}    s_{k-1}s_{k-2}
        \end{cases}$ & $\begin{array}{l}(\lambda_1+\ldots+\lambda_k+\lambda_{k-2}+\ldots+\lambda_2+(2k-4),\lambda_3,\ldots,\lambda_{k-2},\lambda_k,\lambda_{k-1})\\
(\lambda_2+\ldots+\lambda_k+\lambda_{k-2}+\ldots+\lambda_3+(2k-6),\lambda_3+\lambda_2+\lambda_1+2,\lambda_4,\ldots,\lambda_{k-2},\lambda_k,\lambda_{k-1})\\
\quad\vdots\\
(\lambda_k,\lambda_1,\ldots,\lambda_{k-4},\lambda_{k-3}+\ldots+\lambda_k+\lambda_{k-2},\lambda_{k-3})\\
\end{array}$\\
$n-1$ & $\begin{cases}
         s_2s_3\cdots s_ks_{k-2}\cdots s_2s_1\\s_2s_1s_3s_4\cdots s_ks_{k-2}\cdots s_3s_2\\
\vdots\\
	 s_2s_1s_3s_2s_4\cdots 
	 s_{k-2}s_{k-3}     s_{k-1}s_{k-2}s_k
        \end{cases}$ & $\begin{array}{l}(\lambda_1+\ldots+\lambda_k+\lambda_{k-2}+\ldots+\lambda_2+(2k-4),\lambda_3,\ldots,\lambda_{k-2},\lambda_k,\lambda_{k-1})\\
(\lambda_2+\ldots+\lambda_k+\lambda_{k-2}+\ldots+\lambda_3+(2k-6),\lambda_3+\lambda_2+\lambda_1+2,\lambda_4,\ldots,\lambda_{k-2},\lambda_k,\lambda_{k-1})\\
\quad\vdots\\
(\lambda_k,\lambda_1,\ldots,\lambda_{k-4},\lambda_{k-3}+\ldots+\lambda_k+\lambda_{k-2},\lambda_{k-3})\\
\end{array}$\\
\vdots & \;\vdots & \quad\vdots\\
$3k-6$ & $\begin{cases}
           s_2s_3\cdots s_ks_{k-2}\cdots s_2s_1s_2 \cdots s_{k-2}\\s_2s_1s_3s_4\cdots s_ks_{k-2}\cdots s_4s_3s_2s_3\cdots s_{k}\\
\textrm{etc.}
          \end{cases}$ & $\begin{array}{l}(\lambda_1+\ldots+\lambda_k+(k-1),\lambda_2,\ldots,\lambda_{k-3},\lambda_k+\lambda_{k-2}+1,\lambda_{k-1}+\lambda_{k-2}+1)\\(\lambda_2+\ldots+\lambda_{k-1}+(k-3),\lambda_2+\lambda_1+1,\lambda_3,\ldots,\lambda_{k-2},\lambda_k+\lambda_{k-1}+\lambda_{k-2}+2)\\
\textrm{etc.}\end{array}$\\
\vdots & \;\vdots & \quad\vdots\\
$2n-4$ & $s_2\cdots s_ks_{k-2}\cdots s_2s_1s_2 \cdots s_{k-2}s_k\cdots s_3$ & $(\lambda_1+\lambda_2+1,\lambda_3+\lambda_2+1,\lambda_4,\ldots,\lambda_k)$\\
$2n-3$ & $s_2\cdots s_ks_{k-2}\cdots s_2s_1s_2 \cdots s_{k-2}s_k\cdots s_2$ & $(\lambda_1,\lambda_3,\ldots,\lambda_{k})$

\end{tabular}
\end{table}

\normalsize
\end{landscape}

\subsection{Cohomological representations}\label{cohrep}

The determination of data occurring as type $(\pi,w)$ now proceeds to investigating the irreducible cuspidal automorphic representations $\pi$ of $\qL_{\qP}(\ad)$ with non-zero cohomology at the archimedean place, \textit{cohomological representations} for short. 

These are described in terms of the Vogan-Zuckerman classification via its Harish-Chandra modules $A_{\mathfrak{q}}(\lambda)$ (cf. \cite{VZ}), which in turn are parametrised by pairs $(\mathfrak{q},\lambda)$, where $\mathfrak{q}$ is a $\theta$-stable parabolic subalgebra of $\mathfrak{m}=\mathfrak{m}_{\qP}$ and $\lambda$ is an admissible one dimensional representation of the Levi factor $\mathfrak{l}$ of $\mathfrak{q}$. (Throughout the present section this is the only meaning of $\lambda$, and this meaning of $\lambda$ is confined to that one section).
The constructive part of the Vogan-Zuckerman classification states that for $(\mathfrak{q},\lambda)$  there exists a unique irreducible $(\mathfrak{m},K_M^+)$-module $A_{\mathfrak{q}}(\lambda)$ with non-zero cohomology. 

Since the isomorphism class of $A_{\mathfrak{q}}(\lambda)$ is determined by its minimal $K_M^+$-type,
it will suffice to give a complete set of representatives of pairs $(\mathfrak{q},\lambda)$ up to $K_M^+$-conjugacy realising all possible minimal $K_M^+$-types. 
With notation as in section \ref{pre} the minimal $K_M^+$-type of $A_{\mathfrak{q}}(\lambda)$ is equal to $\lambda_{|\mathfrak{t}}+2\rho(\mathfrak{u}\cap\mathfrak{p})$.
So actually, it will suffice to give a complete set of representatives of pairs $(\mathfrak{q},\lambda)$  realising all possible sets $\Phi(\mathfrak{u}\cap\mathfrak{p})$. Also, instead of $\mathfrak{q}\subset\mathfrak{m}$ its Levi subgroup $L=N_M(\mathfrak{q})$, the normaliser of $\mathfrak{q}$ in $M$, will be given.

\subsubsection{The first standard maximal parabolic $\q$-subgroup}
  
Let $\qP=\qP_1$ and drop the index for all objects associated to $\qP$. In particular, write $\qL=\mathbf{SO}(F_{n-1,1})\times \mathbf{GL}_1$ for the Levi component and $M\simeq SO(n-1,1)$ for the group of real points of the semisimple part therein. Recall that the Cartan decomposition corresponding to $K_M$ is written $\mathfrak{m}_0=\mathfrak{k}_0\oplus\mathfrak{p}_0$. Also, notice that an element $X\in\mathfrak{m}_0$ looks like
$$X=\left(\begin{array}{ccc}
a&u&0\\
x&b&-{}^tu\\
0&-{}^tx&-a
\end{array}\right)$$
with $a\in\rea$, $x$ and ${}^tu$ vectors in $\rea^{n-2}$ and $b$ skew-symmetric with respect to 
$I_{n-2}$. 

\subsection*{The odd case}
Assume $2\nmid n$ and set $2(k-1)=n-1$. Then $\mathfrak{m}$ is of type $B_{k-1}$,  
$\textrm{rk}\mathfrak{m}=\textrm{rk}\mathfrak{k}$, and $\mathfrak{b}_0^c=\{H\in\mathfrak{k}_0 \vert b \textrm{ block diagonal of } 2\times 2 \textrm{ blocks}, u=(0,\ldots,0,u_{n-2})\}$ is a maximally compact $\theta$-stable Cartan subalgebra. The linear functional $\varepsilon_2$ now extracts $\sqrt{-2}u_{n-2}$.
Since it is trivial that $\varepsilon_i\in\Phi(\mathfrak{u}\cap\mathfrak{p})$ implies 
$\varepsilon_{i-1}\in\Phi(\mathfrak{u}\cap\mathfrak{p})$ for $3\leq i\leq k$, the possibilities for 
$\Phi(\mathfrak{u}\cap\mathfrak{p})$ are $$\varnothing, \{\varepsilon_2\},\{\varepsilon_2,\varepsilon_3\}\ldots,\{\varepsilon_2,\ldots,\varepsilon_{k-1},\pm\varepsilon_k\}.$$
A complete set of representatives of $\theta$-stable parabolic subalgebras realising all possible 
sets $\Phi(\mathfrak{u}\cap\mathfrak{p})$ is given by any collection of $\theta$-stable parabolic subalgebras $\mathfrak{q}_0,\ldots,\mathfrak{q}_{k-2},\mathfrak{q}_{k-1}^{\pm}$ containing $\mathfrak{b}^c$ and with Levi subgroups isomorphic to the reductive groups listed in table $9$.

\small

\begin{table}[ht]
\caption{Levi subgroups}

\begin{tabular}{c c }
$i$ &  $L_i$\\

\hline

$0$ &  $SO(2(k-1),1)^+$\\
$1$ & $SO(2)\times SO(2(k-2),1)^+$\\
$\vdots$  & $\vdots$\\
$i$ &$SO(2)^i\times SO(2(k-i-1),1)^+$\\
$\vdots$  & $\vdots$\\
$k-1$ & $SO(2)^{k-1}$
\end{tabular}

\end{table}

\normalsize

Let $F=F_{\mu}$ be an irreducible finite dimensional representation of $\mathbf{SO}(n,\co)$ with highest weight $\mu$. For $0\leq i\leq k-2$ there is exactly one irreducible unitary representation $A_i(\lambda)$ of $SO(n-1,1)^+$ with non-zero cohomology twisted by $F$, the one with $\lambda_{|\mathfrak{b}^c}=\mu$, and  
$$H^q(\mathfrak{m}, K^+_{M}; A_i(\lambda)\otimes F_{\mu})=\left\{ \begin{array}{ll}\co& q=i,n-1-i, 
\lambda_{|\mathfrak{b}^c}=\mu\\0 & \textrm{otherwise}.\end{array}\right.$$

For $i=k-1$ on the other hand there are two irreducible unitary representations $A_{k-1}^{\pm}(\lambda)$ of $SO(n-1,1)^+$ with non-zero cohomology twisted by $F$ corresponding to 
whether $\Phi(\mathfrak{u}\cap\mathfrak{p})$ contains $\varepsilon_k$ or $-\varepsilon_k$.

\begin{rem}
Since the Levi subgroup $L_i$ is compact for $i=k-1$, the two representations $A_{k-1}^{\pm}(\lambda)$ belong to the discrete series of $M^+$.
\end{rem}

If the highest weight $\mu$ of $F$ is assumed to be regular the necessary condition $\lambda_{|\mathfrak{b}^c}=\mu$ can be met only in case of the two representations $A_{k-1}^{\pm}(\lambda)$.
As for a representation of $M=SO(n-1,1)$ consider the induced representation 
$\textrm{Ind}_{M^+}^M(A_{k-1}^+)=\textrm{Ind}_{M^+}^M(A_{k-1}^-)$, the restriction of 
which back to $M^+$ gives the direct sum $A_{k-1}^{+}(\lambda)\oplus A_{k-1}^{-}(\lambda)$. The final result is summarised in
\begin{prop}\label{car1odd}
Let $F$ be an irreducible finite dimensional representation of $\mathbf{SO}(n,\co)$ of regular highest weight $\mu$, also denoted by $F_{\mu}$. The irreducible quotient $A_{k-1}(\lambda)$ of the induced representation  $\textrm{Ind}_{M^+}^M(A_{k-1}^+(\lambda))=\textrm{Ind}_{M^+}^M(A_{k-1}^-(\lambda))$, where $\lambda=\mu$, then is the unique irreducible unitary representation of $SO(n-1,1)$ with non-zero cohomology twisted by $F$ and
$$H^q(\mathfrak{m}, K_M; A_{k-1}(\lambda)\otimes F_{\mu})=\left\{ \begin{array}{ll}\co& q=k-1, \lambda=\mu\\0 & \emph{\textrm{otherwise}}.\end{array}\right.$$
\end{prop}

\subsection*{The even case}

Let now $2\mid n$ and set $2(k-1)=n$. Then $\mathfrak{m}$ is of type $D_{k-1}$ and $\textrm{rk}\mathfrak{m}=\textrm{rk}\mathfrak{k}+1$.
The real form $\mathfrak{b}_0$ of the Cartan subalgebra of $\mathfrak{m}$ chosen in \ref{Levis} is both maximally compact and maximally split. Here, the linear functionals introduced above are ordered 
as $\varepsilon_3\geq\varepsilon_4\geq\ldots\geq\varepsilon_k\geq\varepsilon_2$. (The reason for the reordering is that the positive system 
to be defined on $\Phi(\mathfrak{m},\mathfrak{b})$ needs to be $\theta$-stable).
In analogy with the odd case, the possibilities for $\Phi(\mathfrak{u}\cap\mathfrak{p})$ are
$$\varnothing, \{\varepsilon_3\},\ldots,\{\varepsilon_3,\ldots,\varepsilon_{k}\}.$$
A complete set of representatives of $\theta$-stable parabolic subalgebras realising all possible 
sets $\Phi(\mathfrak{u}\cap\mathfrak{p})$ is given by any collection of $\theta$-stable parabolic subalgebras $\mathfrak{q}_0,\ldots,\mathfrak{q}_{k-2}$ containing $\mathfrak{b}$ and with Levi subgroups isomorphic to the reductive groups listed in table $10$. 

\small

\begin{table}[ht]
\caption{Levi subgroups}

\begin{tabular}{c c }
$i$ &   $L_i$\\

\hline

$0$ &   $SO(2k-3,1)^+$\\
$1$ & $SO(2)\times SO(2k-5,1)^+$\\
$\vdots$  & $\vdots$\\
$i$ & $SO(2)^i\times SO(2k-2i-3,1)^+$\\
$\vdots$  & $\vdots$\\
$k-3$ & $SO(2)^{k-3}\times SO(3,1)^+$\\
$k-2$ & $SO(2)^{k-2}\times SO(1,1)^+$
\end{tabular}

\end{table}

\normalsize

\begin{rem}
The groups displayed do not exhaust the Levi subgroups of non-isomorphic $\theta$-stable parabolic 
subalgebras of $\mathfrak{m}$; e.g. there is always one being isomorphic to $SO(2)^{k-2}\times SO(2,1)^+$.
Still, they parametrise all non-isomorphic irreducible $(\mathfrak{m}, K^+_{M})$-modules
with non-zero cohomology determined by the Vogan-Zuckerman classification. 
\end{rem}

Let $F=F_{\mu}$ be an irreducible finite dimensional representation of $\mathbf{SO}(n,\co)$ with highest weight 
$\mu$. For $0\leq i\leq k-2$ there is exactly one irreducible unitary representation $A_i(\lambda)$ of $SO(n-1,1)^+$ with non-zero cohomology twisted by $F$, the one where $\lambda_{|\mathfrak{b}}=\mu$, and
$$H^q(\mathfrak{m}, K^+_{M}; A_i(\lambda)\otimes F_{\mu})=\left\{ \begin{array}{ll}\co& q=i,n-1-i, \begin{array}{l}\lambda_{|\mathfrak{b}}=\mu \textrm{ if } 2\mid k-1\\\lambda_{|\mathfrak{b}}=\eta(\mu) \textrm{ if } 2\nmid k-1\end{array}\\0 & \textrm{otherwise}.\end{array}\right.$$
(Here, $\eta$ denotes the Dynkin diagram automorphism induced by the unique element of maximal length in $W$, interchanging $\alpha_{k-1}$ and $\alpha_k$).
\begin{rem}
The representation $A_{k-2}(\lambda)$ is tempered, unitarily induced.
\end{rem}

Assume the highest weight $\mu$ of $F$ now to be regular and consider the induced representation $A_{k-2}(\lambda)=\textrm{Ind}_{M^+}^M(A_{k-2}(\lambda))$ of $M$. The final result is contained in 
\begin{prop}\label{car1even}
Let $F$ be an irreducible finite dimensional representation of $\mathbf{SO}(n,\co)$ of regular highest weight $\mu=(\mu_2,\ldots,\mu_k)$, also denoted by $F_{\mu}$.
If $\mu_{k}=\mu_{k-1}$, the irreducible quotient $A_{k-2}(\lambda)$ of $\textrm{Ind}_{M^+}^M(A_{k-2}(\lambda))$, where $\lambda=\mu$ is the unique irreducible unitary representation of $SO(n-1,1)$ with non-zero cohomology twisted by $F$ and
$$H^q(\mathfrak{m}, K_M; A_{k-2}(\lambda)\otimes F_{\mu})=\left\{ \begin{array}{ll}\co& q=k-2,k-1, \lambda=\mu \\0 & \emph{\textrm{otherwise}}.\end{array}\right.$$
If $\mu_{k}\neq\mu_{k-1}$, there is no irreducible unitary representation of $SO(n-1,1)^+$ with non-zero cohomology twisted by $F$.                                                                                 
 \end{prop}

\begin{rem}
Rewriting $\mu=(\mu_2,\ldots,\mu_k)$ in coordinates for the linear functionals $\varepsilon_i$ get
\small
$$\mu=(\mu_2+\ldots+\mu_{k-2}+\frac{1}{2}(\mu_{k-1}+\mu_k))\varepsilon_2+\ldots+\frac{1}{2}(\mu_{k-1}+\mu_k)
\varepsilon_{k-1}+\frac{1}{2}(-\mu_{k-1}+\mu_k)\varepsilon_k.$$
\normalsize
Doing the same thing for an admissible character $\lambda$ of $L_i$, $i\in\{1,\ldots,k-2\}$
yields $$\lambda=n_2\varepsilon_2+\ldots+n_{k-1}\varepsilon_{k-1}$$
with $n_i\in\z$, since $\lambda$ is assumed to be unitary. 
Given that $\mu$ is regular the requirement $\lambda|_{\mathfrak{b}}=\mu$ can be met only if on the one hand $\lambda$ is
admissible with respect to $\Phi^+(\mathfrak{m},\mathfrak{t})$, that is $n_i\neq0$ for $2\leq i\leq k-1$,
and on the other $-\mu_{k-1}+\mu_k=0$.
\end{rem}

\subsubsection{The second standard maximal parabolic $\q$-subgroup}
 
Let $\qP=\qP_2$ and drop the index for all objects associated to $\qP$. In particular, $\qL=\mathbf{SO}(F_{n-2})\times \mathbf{GL}(2)$ shall denote the Levi component of $\qP$ and $M=SL_2^{\pm}(\rea)\times SO(n-2)$ the group of real points of the semisimple part therein.
 
Here, the description of cohomological representations does not depend on the parity of $n$. Thus, the results are stated in a way appealing to both the odd and even case, though still in terms of pairs
$(\mathfrak{q},\lambda)$.

\subsection*{The odd case and even case}
In both cases there are just two Levi subgroups up to isomorphism, $L_0=SL_2(\rea)\times SO(n-2)$ and $L_{k-1}=SO(2)^{k-1}$. If $F$ is an irreducible finite dimensional representation of 
$\mathbf{SL}_2(\co)\times\mathbf{SO}(n-2,\co)$ of regular highest weight $\mu$, in case of $L_{k-1}$ there are two irreducible unitary representations $A_{k-1}^+(\lambda)$ and $A_{k-1}^-(\lambda)$ of $M^+$ with non-zero cohomology twisted by $F$. They owe to the discrete series of $SL_2(\rea)$. The  regularity of $\mu$ entails that there is none for $L_0$. The result for $M$ is stated in 

\begin{prop}\label{car2odd}
Let $F=F_{\mu}$ be an irreducible finite dimensional representation of 
$\mathbf{SL}_2(\co)\times\mathbf{SO}(n-2,\co)$ of regular highest weight $\mu$. There is a unique irreducible unitary representation $A_{k-1}(\lambda)$ of $SL_2^{\pm}(\rea)\times SO(n-2)$ with non-zero cohomology twisted by $F$, necessarily with $\lambda=\mu$, and
$$H^q(\mathfrak{m}, K_M; A_{k-1}(\lambda)\otimes F_{\mu})=\left\{ \begin{array}{ll}\co& q=1, \lambda=\mu\\0 & \emph{\textrm{otherwise}}.\end{array}\right.$$
This representation is distinguished by the property that its restriction to $M^+$ yields the direct sum  $A_{k-1}^+(\lambda)\oplus A_{k-1}^-(\lambda)$ of the two discrete series representations $A_{k-1}^+(\lambda)$ and $A_{k-1}^-(\lambda)$ of $M^+$.
\end{prop}

\subsection{Classes of type $(\pi,w)$}

The results obtained earlier in this section are combined and the classes of type $(\pi,w)$ for the two standard maximal parabolic $\q$-subgroups of $\G$ are listed.

\subsubsection*{The first standard maximal parabolic $\q$-subgroup}

Let $\qP=\qP_1$ and drop the index for all objects associated to and distinct from $\qP$.

\subsection*{The odd case}

Let $F_{\mu_w}$ be the irreducible finite dimensional representation of $\M(\co)\simeq\mathbf{SO}(n,\co)$ of 
highest weight ${\mu_w}_{|\mathfrak{b}}$ from Kostant's theorem (cf. \ref{type}). Moreover, let $\delta$ be the Cayley transform mapping $\mathfrak{b}$ to $\mathfrak{b}^c$, the complexification of the maximally compact $\theta$-stable Cartan subalgebra. 
According to Proposition \ref{car1odd} the representation $A_{k-1}(\lambda)$ where $\lambda=\delta({\mu_w}_{|\mathfrak{b}})$ is the unique irreducible unitary representation of $M\simeq SO(n-1,1)$ with non-zero cohomology twisted by $F_{\mu_w}$. 

Henceforth, notation is switched from the Harish-Chandra module to the unitary representation in order 
to introduce the abbreviation $\pi_{k-1}(\mu)=\pi(A_{k-1}(\delta({\mu_w}_{|\mathfrak{b}})))$. 
\begin{quotation}
The possible types of cuspidal classes for $\qP_1$ in the odd case then are $(\pi_{k-1}(\mu),w)$ for $w$ ranging through $W^{\qP_1}$. 
\end{quotation}

\subsection*{The even case}

Let $F_{\mu_w}$ be the irreducible finite dimensional representation of $\M(\co)\simeq\mathbf{SO}(n,\co)$ 
of highest weight ${\mu_w}_{|\mathfrak{b}}=(\mu_2,\ldots,\mu_k)$ from Kostant's theorem and suppose $\mu_{k-1}=\mu_k$. In view of table $6$ from \ref{Lmodules} this constraint excludes the two representations
$F_{\mu_w}$ with $l(w)=k-1$ from the considerations to follow and amounts to requiring $\lambda_{k-1}=\lambda_k$ 
for the highest weight $\lambda=(\lambda_1,\dots,\lambda_k)$ of $(\tau,E)$ in case of the remaining ones. 

Here, the $\theta$-stable Cartan subalgebra $\mathfrak{b}_0$ is itself maximally compact. Therefore, according to Proposition \ref{car1even} the representation $A_{k-2}(\lambda)$ where $\lambda={\mu_w}_{|\mathfrak{b}}$ is the unique irreducible unitary representation of $M\simeq SO(n-1,1)$ with non-zero cohomology twisted by $F_{\mu_w}$. 
Again using the abbreviation $\pi_{k-2}(\mu)=\pi(A_{k-2}({\mu_w}_{|\mathfrak{b}}))$
\begin{quotation} the possible types of cuspidal classes for $\qP_1$ in the even case are $(\pi_{k-2}(\mu),w)$ for $w$ ranging through $W^{\qP_1}$. 
\end{quotation}

\subsubsection*{The second standard maximal parabolic $\q$-subgroup}

Let $\qP=\qP_2$ and drop the index for all objects associated to and distinct from $\qP$. 
Since the subsequent discussion is formally identical for both cases it is given but once.

Let $F_{\mu_w}$ now be the irreducible finite dimensional representation of $\M(\co)\simeq\mathbf{SL}_2(\co)\times\mathbf{SO}(n-2,\co)$ of highest weight ${\mu_w}_{|\mathfrak{b}}$. 
Moreover, let $\delta$ be the Cayley transform mapping $\mathfrak{b}$ to $\mathfrak{b}^c$, the complexification 
of the maximally compact $\theta$-stable Cartan subalgebra. 
According to Proposition \ref{car2odd} in either case there is exactly one irreducible unitary representation $A_{k-1}(\lambda)$ of $M\simeq SL_2^{\pm}(\rea)\times SO(n-2)$ with non-zero cohomology twisted by $F_{\mu_w}$, the
one with $\lambda=\delta({\mu_w}_{|\mathfrak{b}})$. 
In the notation $\pi_{k-1}(\mu)=A_{k-1}(\delta({\mu_w}_{|\mathfrak{b}}))$ 
\begin{quotation}
the possible types of cuspidal classes for $\qP_2$ are $(\pi_{k-1}(\mu),w)$ for $w$ ranging through $W^{\qP_2}$.
\end{quotation}

\section{Regular Eisenstein cohomology classes}\label{construction}

This section provides the construction of regular Eisenstein cohomology classes for the two standard maximal parabolic $\q$-subgroups of $\G$. The theorems appealing to the generic case (cf. section \ref{regular}) are applied to obtain the final results.
As a delineation of possible degrees in which those regular Eisenstein cohomology classes occur  the general vanishing theorems are employed first.

Suppose now that the highest weight $\lambda_{\tau}=(\lambda_1,\dots,\lambda_k)$ of the irreducible finite dimensional complex rational representation $(\tau,E)$ be regular.

\subsection{The range of cohomology}

The values for the two quantities $l_0$ and $q_0$ associated to $\G$ are computed to be $l_0(\G(\rea))=0$ and $q_0(\G(\rea))=n$. As a consequence,  the following restrictions on the cohomology of $\G$ can be inferred:

\begin{enumerate}
\item the automorphic cohomology of $\G$ vanishes below degree $n$, i.e. $H^q(\mathfrak{g},K;\mathcal{A}_E\otimes E)=0$ for $q<n$. For $q\geq n$ each summand $H^q(\mathfrak{g},K;\mathcal{A}_{E,\{\qP\}}\otimes E)$, $\{\qP\}\neq\{\G\}$ is generated by \textit{regular} Eisenstein cohomology classes. 
\item the cuspidal cohomology of $\G$ occurs only in degree $n$, i.e. $H^q(\mathfrak{g},K;{}^{\circ}\mathcal{A}_E\otimes E)=0$ for $q\neq n$.
\end{enumerate}
The virtual cohomological dimension is found to be $\textrm{vcd}(\G)=2n-2$. Hence, as a matter of 
\begin{fact}
In the generic case the automorphic cohomology $H^q(\mathfrak{g},K;\mathcal{A}_E\otimes E)$ of $\G$ vanishes for $q\not\in[n,2n-2]$.
\end{fact}

\subsection{Regular Eisenstein cohomology classes}

Let $\qP$ be a standard maximal parabolic $\q$-subgroup, $\pi$ a cohomological representation of 
$\qL_{\qP}(\ad)$ and $f\in W_{\qP,\tilde \pi}$. The Eisenstein series $E^{\G}_{\qP}(f,\lambda)$ is 
evaluated at the real points $\lambda_w=-w(\lambda_{\tau}+\rho)|_{\mathfrak{a}_P}$ , $w\in W^{\qP}$ 
to obtain non-trivial cohomology classes.

\subsubsection*{The first standard maximal parabolic $\q$-subgroup} 
For the sake of transparency the index is maintained when stating the final results now.

\subsection*{The odd case}
Let $2\nmid n$ and $2k=n+1$. The real parameters $\lambda_w$ for $w\in W^{\qP_1}$ are given in
their normalisation with respect to $\rho_1$ in table $11$.

\scriptsize

\begin{table}[ht]
\caption{$\lambda_w=a_1\rho_1$}

\begin{tabular}{c l l}
$l(w)$ & $w$ & $a_1$\\
\hline

$0$ & $1$ & $-(\lambda_1+\ldots+\lambda_k+\ldots+\lambda_1+n)\cdot \frac{1}{n}$\\
$1$ & $s_1$ & $-(\lambda_2+\ldots+\lambda_k+\ldots+\lambda_2+(n-2))\cdot \frac{1}{n}$\\
\vdots & \;\vdots & \quad\vdots\\
$k-1$ & $s_1\dots s_{k-1}$ &  $-(\lambda_k+1)\cdot \frac{1}{n}$\\
$k$ & $s_1\dots s_{k-1}s_k$ &  \: $(\lambda_k+1)\cdot \frac{1}{n}$\\
\vdots & \;\vdots & \quad\vdots\\
$n-1$ & $s_1\dots s_{k-1}s_k\dots s_2$ & \: $(\lambda_2+\ldots+\lambda_k+\ldots+\lambda_2+(n-2))\cdot \frac{1}{n}$\\
$n$ & $s_1\dots s_{k-1}s_k\dots s_2s_1$ & \: $(\lambda_1+\ldots+\lambda_k+\ldots+\lambda_1+n)\cdot \frac{1}{n}$
\end{tabular}

\end{table}
\normalsize

The final result in the present case may now be stated as the

\begin{thm}\label{myfirst}
Let $\G$ be the $\q$-rational form of $\q$-rank $2$ of $SO(n,2)$ introduced in section \ref{q-form}. Suppose the  highest weight $\lambda_{\tau}=(\lambda_1,\ldots,\lambda_k)$ of $(\tau,E)$ to be regular. Let $\qP_1$ be the standard maximal parabolic $\q$-subgroup of $\G$ with
Levi component $\qL_1\simeq\mathbf{SO}(F_{n-1,1})\times\mathbf{GL}_1$. In the \emph{odd} case the $\G(\af)$-module $H^q(\mathfrak{g},K;\mathcal{A}_{E,\{\qP_1\}}\otimes E)$ 
\begin{itemize}
\item[$\circ$] is generated in degree $q=k-1+l$ for $k\leq l\leq n$ by the regular Eisenstein cohomology classes $[E^{\G}_{\qP_1}(f,\lambda_w)]$ associated to cuspidal classes of the single type $(\pi_{k-1}(\mu),w)$ with $l(w)=l$
\item[$\circ$] vanishes otherwise, that is for $q\not\in[n,\frac{3n-1}{2}]$.
\end{itemize}
\end{thm}

\begin{rem}
 Notice that $4n-4-(3n-1)=n-3$, which for $n\geq5$ implies that the summand
 $H^q(\mathfrak{g},K;\mathcal{A}_{E,\{\qP_1\}}\otimes E)$ doesn't ever contribute to the 
 automorphic cohomology of $\G$ in degree $q=\textrm{vcd}(\G)$.  
\end{rem}

\begin{proof}
Take a cuspidal class of type $(\pi_{k-1}(\mu),w)$, the degree of which is $k-1$, and consider the 
associated Eisenstein series $E^{\G}_{\qP_1}(f,\lambda)$ for $f\in W_{\qP_1,\tilde\pi_{k-1}(\mu)}$. According to Theorem \ref{mainthm}, section \ref{regular}, if 
$E^{\G}_{\qP_1}(f,\lambda)$ is holomorphic at $\lambda_w$, the automorphic form $E^{\G}_{\qP_1}(f,\lambda_w)$ will represent a non-trivial 
Eisenstein cohomology class of degree $q=k-1+l(w)$. 

By the vanishing of the automorphic cohomology of $\G$ in degrees less than $n$ only the classes of type $(\pi_{k-1}(\mu),w)$ for those $w$ have to be taken into account for which $k-1+l(w)\geq n$, that is for which $l(w)\geq k$. Since $\qP_1$ is maximal and $k\geq\frac{1}{2}\dim\mathfrak{n}_1$, Theorem 
\ref{holomorphy}, section \ref{regular}, ensures that $E^{\G}_{\qP_1}(f,\lambda)$ is holomorphic at $\lambda_w$.

Finally, $k-1+n=\frac{3n-1}{2}$ computes the upper bound for the non-vanishing.
\end{proof}

\subsection*{The even case}
 
Let $2\mid n$ and $2k=n+2$. The real parameters $\lambda_w$ for $w\in W^{\qP_1}$ are given in
their normalisation with respect to $\rho_1$ in table $12$.

\scriptsize

\begin{table}[ht]
\caption{$\lambda_w=a_1\rho_1$}

\begin{tabular}{c l l}
$l(w)$ & $w$ & $a_1$\\
\hline

$0$ & $1$ & $-(\lambda_1+\ldots+\lambda_k+\lambda_{k-2}+\ldots+\lambda_1+n)\cdot\frac{1}{n}$\\
$1$ & $s_1$ & $-(\lambda_2+\ldots+\lambda_k+\lambda_{k-2}+\ldots+\lambda_2+(n-2))\cdot\frac{1}{n}$\\
\vdots & \;\vdots & \quad \vdots  \\
$k-2$ & $s_1\dots s_{k-2}$ &  $-(\lambda_{k-1}+\lambda_k+2)\cdot\frac{1}{n}$\\
$k-1$ & $\begin{cases}s_1\dots s_{k-2}s_k\\s_1\dots s_{k-2}s_{k-1}\end{cases}$ &  $\begin{array}{l}-(\lambda_{k-1}-\lambda_k)\cdot\frac{1}{n}\\(\lambda_{k-1}-\lambda_k)\cdot\frac{1}{n}                                                                                   \end{array}$\\
$k$ & $s_1\dots s_{k-2}s_{k-1}s_k$ & \: $(\lambda_{k-1}+\lambda_k+2)\cdot\frac{1}{n}$\\
\vdots & \;\vdots & \quad \vdots \\
$n-1$ & $s_1\dots s_{k-1}s_ks_{k-2}\dots s_2$ & \: $(\lambda_2+\ldots+\lambda_k+\lambda_{k-2}+\ldots+\lambda_2+(n-2))\cdot\frac{1}{n}$\\
$n$ & $s_1\dots s_{k-1}s_ks_{k-2}\dots s_2s_1$ & \: $(\lambda_1+\ldots+\lambda_k+\lambda_{k-2}+\ldots+\lambda_1+n)\cdot\frac{1}{n}$
\end{tabular}

\end{table}
\normalsize

The final result about the summand in the automorphic cohomology of $\G$ indexed by $\qP_1$ in the even case is

\begin{thm}
Let $\G$ be the $\q$-rational form of $\q$-rank $2$ of $SO(n,2)$ introduced in section \ref{q-form}. Suppose the  highest weight $\lambda_{\tau}=(\lambda_1,\ldots,\lambda_k)$ of $(\tau,E)$ to be regular and such that $\lambda_{k-1}=\lambda_k$. Let $\qP_1$ be the standard maximal parabolic $\q$-subgroup of $\G$ with
Levi component $\qL_1\simeq\mathbf{SO}(F_{n-1,1})\times\mathbf{GL}_1$. In the \emph{even} case the $\G(\af)$-module $H^q(\mathfrak{g},K;\mathcal{A}_{E,\{\qP_1\}}\otimes E)$ 
\begin{itemize}
\item[$\circ$] is generated in degrees $q\in[n,\frac{3n}{2}]$
by the regular Eisenstein cohomology 
classes $[E^{\G}_{\qP_1}(f,\lambda_w)]$ associated to the cuspidal classes of type $(\pi_{k-2}(\mu),w)$ of degree $d=k-2,k-1$ with $l(w)\in\{k,\ldots,n\}$ such that $q=d+l(w)$
\item[$\circ$] vanishes otherwise, that is for $q\not\in[n,\frac{3n}{2}]$.
\end{itemize}
\end{thm}

\begin{rem}
 (i) As $n\geq6$ the difference $2\cdot\textrm{vcd}(\G)-3n=n-4$ shows that the summand $H^q(\mathfrak{g},K;\mathcal{A}_{E,\{\qP_1\}}\otimes E)$  doesn't ever contribute in degree $q=\textrm{vcd}(\G)$.
 
(ii) The proof is analogous to the odd case and as $k-1+n=\frac{n}{2}+n$ also the vanishing statement 
is seen to be true.
\end{rem}

\subsubsection*{The second standard maximal parabolic $\q$-subgroup}

As before, the index is maintained for stating the final results.

The real parameters $\lambda_w$ for $w\in W^{\qP_2}$ are given in their normalisation with respect to $\rho_2$ in tables $13$ and $14$ for the odd ($2k=n+1$) and the even case ($2k=n+2$) respectively.

\scriptsize

\begin{table}[ht]\label{lambdaw3}
\caption{$\lambda_w=a_2\rho_2$}
\begin{tabular}{c l l}
$l(w)$ & $w$ & $a_2$\\
\hline

$0$ & $1$ & $-(\lambda_1+\ldots+\lambda_{k}+\ldots+\lambda_2+(n-1))\cdot\frac{1}{n-1}$\\
$1$ & $s_2$ & $-(\lambda_1+\ldots+\lambda_k+\ldots\lambda_3+(n-2))\cdot \frac{1}{n-1}$\\
\vdots & \;\vdots & \quad\vdots\\
$k-2$ & $\begin{cases}s_2\cdots s_{k-1}\\s_2s_1s_3s_4\cdots s_{k-2}\\
\textrm{etc.} 
          \end{cases}$
      & $\begin{array}{l}
          -(\lambda_1+\ldots+\lambda_{k}+(n-k+1))\cdot \frac{1}{n-1}\\
	  -(\lambda_2+\ldots+\lambda_{k}+\lambda_{k-1}+(n-k+1))\cdot \frac{1}{n-1}\\
	  \textrm{etc.}
         \end{array}$\\
\vdots & \;\vdots & \quad\vdots\\
$n-2$ & $\begin{cases}s_2\cdots s_k\cdots s_2\\s_2s_1s_3s_4\cdots s_k\cdots s_4s_3\\
\vdots\\s_2s_1s_3s_2s_4s_3\cdots s_{k-1}s_{k-2}s_k          
         \end{cases}$ &
          $\begin{array}{l}
          -(\lambda_1+1)\cdot \frac{1}{n-1}\\
	  -(\lambda_2+1)\cdot \frac{1}{n-1}\\
\vdots\\-(\lambda_{k-1}+1)\cdot \frac{1}{n-1}
	           \end{array}$\\
$n-1$ & $\begin{cases}s_2\cdots s_k\cdots s_2s_1\\s_2s_1s_3s_4\cdots s_k\cdots s_4s_3s_2\\
\vdots\\s_2s_1s_3s_2s_4s_3\cdots s_{k-2}s_ks_{k-1}       
         \end{cases}$ &
          $\begin{array}{l}
          (\lambda_1+1)\cdot \frac{1}{n-1}\\
	  (\lambda_2+1)\cdot \frac{1}{n-1}\\
\vdots\\(\lambda_{k-1}+1)\cdot \frac{1}{n-1}
	           \end{array}$\\
\vdots & \;\vdots & \quad\vdots\\
$3k-3$ & $\begin{cases}s_2\cdots s_k\cdots s_2s_1s_2\cdots s_k\\s_2s_1s_3s_4\cdots s_k\cdots s_2\cdots s_ks_{k-1}\\
\textrm{etc.} 
          \end{cases}$
      & $\begin{array}{l}
          (\lambda_1+\ldots+\lambda_{k}+(n-k+1))\cdot \frac{1}{n-1}\\
	  (\lambda_2+\ldots+\lambda_{k}+\lambda_{k-1}+(n-k+1))\cdot \frac{1}{n-1}\\
	  \textrm{etc.}
         \end{array}$\\
\vdots & \;\vdots & \quad\vdots\\
$2n-4$ & $s_2\dots s_k\dots s_2s_1s_2\dots s_k\dots s_3$ & 
$(\lambda_1+\ldots+\lambda_k+\ldots\lambda_3+(n-2))\cdot \frac{1}{n-1}$\\
$2n-3$ & $s_2\cdots s_k\cdots s_2s_1s_2\cdots s_k\cdots s_2$ & $(\lambda_1+\ldots+\lambda_{k}+\ldots+\lambda_2+(n-1))\cdot\frac{1}{n-1}$
\end{tabular}

\end{table}
\normalsize

\begin{rem}
As a matter of fact, which is also suggested by the table, for $\lambda_w=a_2\rho_2$ the minimum to be attained by $a_2$ as $\lambda$ ranges 
through the set of dominant regular weights depends only on the length of $w$. 
\end{rem}

The final result in both cases -- apparently with a shift in meaning for $k$ -- reads as 
\begin{thm}\label{mythird}
Let $\G$ be the $\q$-rational form of $\q$-rank $2$ of $SO(n,2)$ introduced in section \ref{q-form}.  Suppose the highest weight $\lambda_{\tau}=(\lambda_1,\ldots,\lambda_k)$ of $(\tau,E)$ to be regular. Let $\qP_2$ be the standard maximal parabolic $\q$-subgroup of $\G$ with
Levi component $\qL_2\simeq\mathbf{SO}(F_{n-2})\times \mathbf{GL}_2$. In the \emph{odd} case the 
$\G(\af)$-module $H^q(\mathfrak{g},K;\mathcal{A}_{E,\{\qP_2\}}\otimes E)$ 
\begin{itemize}
\item[$\circ$] is generated in degree $q=1+l$ for $n-1\leq l\leq 2n-3$ by the regular Eisenstein cohomology classes $[E^{\G}_{\qP_2}(f,\lambda_w)]$ associated to cuspidal classes of type $(\pi_{k-1}(\mu),w)$ with $l(w)=l$
\item[$\circ$] vanishes otherwise, that is for $q\not\in[n,2n-2]$.
\end{itemize}
\end{thm}

\begin{rem}
(i) Notice that $H^q(\mathfrak{g},K;\mathcal{A}_{E,\{\qP_2\}}\otimes E)$ contributes in every single possible  degree to the automorphic cohomology of $\G$. 

(ii) The proof follows along the same lines as that of Theorem \ref{myfirst}. Here, $1+2n-3=\textrm{vcd}(\G)$ computes the upper bound for the non-vanishing of $H^q(\mathfrak{g},K;\mathcal{A}_{E,\{\qP_2\}}\otimes E)$.
\end{rem}

\begin{landscape}

\scriptsize
\begin{table}[ht]\label{lambdaw4}
\caption{$\lambda_w=a_2\rho_2$}
\begin{tabular}{c l l}
$l(w)$ & $w$ & $a_2$\\
\hline

$0$ & $1$ & $-(\lambda_1+\ldots+\lambda_k+\lambda_{k-2}+\ldots+\lambda_2+(n-1))\cdot\frac{1}{n-1}$\\
$1$ & $s_2$ & $-(\lambda_1+\ldots+\lambda_k+\lambda_{k-2}+\ldots+\lambda_3+(n-2))\cdot\frac{1}{n-1}$\\
\vdots & \;\vdots & \quad\vdots\\
$k-3$ & $\begin{cases}s_2\cdots s_{k-2}\\s_2s_1s_3s_4\cdots s_{k-3}\\\textrm{etc.}\end{cases}$ & $\begin{array}{l}-(\lambda_1+\ldots+\lambda_k+(n-k+2))
                           \cdot\frac{1}{n-1}\\-(\lambda_2+\ldots+\lambda_k+\lambda_{k-1}+(n-k+2))\cdot\frac{1}{n-1}\\
			   \textrm{etc.}
			   \end{array}$\\
\vdots & \;\vdots & \quad\vdots\\
$n-2$ & $\begin{cases}
         s_2s_3\cdots s_ks_{k-2}\cdots s_3s_2\\s_2s_1s_3s_4\cdots s_ks_{k-2}\cdots s_4s_3\\
\vdots\\
	 s_2s_1s_3s_2s_4\cdots s_{k-2}s_{k-3}  s_{k-1}s_{k-2}
        \end{cases}$ & $\begin{array}{l}-(\lambda_1+1)\cdot\frac{1}{n-1}\\-(\lambda_2+1)\cdot\frac{1}{n-1}\\
\quad\vdots\\
-(\lambda_k+1)\cdot\frac{1}{n-1}
\end{array}$\\
$n-1$ & $\begin{cases}
         s_2s_3\cdots s_ks_{k-2}\cdots s_2s_1\\s_2s_1s_3s_4\cdots s_ks_{k-2}\cdots s_3s_2\\
\vdots\\
	 s_2s_1s_3s_2s_4\cdots s_{k-2}s_{k-3}    s_{k-1}s_{k-2}s_k
        \end{cases}$ & $\begin{array}{l}(\lambda_1+1)\cdot\frac{1}{n-1}\\(\lambda_2+1)\cdot\frac{1}{n-1}\\
\quad\vdots\\
(\lambda_k+1)\cdot\frac{1}{n-1}
\end{array}$\\
\vdots & \;\vdots & \quad\vdots\\
$3k-4$ & $\begin{cases}
           s_2s_3\cdots s_ks_{k-2}\cdots s_2s_1s_2 \cdots s_{k-2}s_ks_{k-1}\\s_2s_1s_3s_4\cdots s_ks_{k-2}\cdots s_4s_3s_2s_3\cdots s_{k-2}s_{k}s_{k-1}s_{k-2}\\
\textrm{etc.}
          \end{cases}$ & $\begin{array}{l}(\lambda_1+\ldots+\lambda_k+(n-k+2))
                           \cdot\frac{1}{n-1}\\(\lambda_2+\ldots+\lambda_k+\lambda_{k-1}+(n-k+2))\cdot\frac{1}{n-1}\\
			   \textrm{etc.}
			   \end{array}$\\
\vdots & \;\vdots & \quad\vdots\\
$2n-4$ & $s_2\cdots s_ks_{k-2}\cdots s_2s_1s_2 \cdots s_{k-2}s_k\cdots s_3$ & $(\lambda_1+\ldots+\lambda_k+\lambda_{k-2}+\ldots+\lambda_3+(n-2))\cdot\frac{1}{n-1}$ \\
$2n-3$ & $s_2\cdots s_ks_{k-2}\cdots s_2s_1s_2 \cdots s_{k-2}s_k\cdots s_2$ & $(\lambda_1+\ldots+\lambda_k+\lambda_{k-2}+\ldots+\lambda_2+(n-1))\cdot\frac{1}{n-1}$

\end{tabular}
\end{table}

\normalsize

\end{landscape}

\begin{appendix}
\section{The set of minimal coset representatives as Hasse diagram}\label{diagram}

Studying reductive algebraic groups $\G$ in general there is the question of how to determine the set $W^{\qP}$ of minimal coset representatives associated to a parabolic $\q$-subgroup $\qP$ by way of a preferably efficient algorithm. Such an algorithm is obtained by elaborating on the results of \cite{Kos}
(cf. \cite{B-E}).

\subsection{The Bruhat order on the Weyl group} Notation is as in section \ref{pre}.
Let $W=W(\G,\T)$ denote the Weyl group for a fixed maximal torus $\T\subset\G$, $\mathbf{B}\subset\G$ be a Borel subgroup containing $\T$ and $\Phi^+$ the corresponding system of positive roots.
Then $W$ can be given the structure of a partially ordered set by introducing the 
\textit{Bruhat order}:
\begin{quote}
set $w'\leq w$ if $w=w''w'$ for some word $w''$ in the alphabet $\{s_{\alpha}\}_{\alpha\in\Phi^+}$ and $l(w')\leq l(w)$,
\end{quote}
where $l(w)$ is the length of $w\in W$.
Depicting $W$ together with the Bruhat order as a directed graph yields what is called a \textit{Hasse diagram}\footnote{See for instance R. Stanley, \emph{Enumerative combinatorics. Vol.1}. Cambridge Studies in Advanced Mathematics. 49. Cambridge: Cambridge University Press. xi, 326 p, 1999.}:
\begin{itemize}
\item[$\circ$] draw a vertex for each element $w\in W$ and
\item[$\circ$] join two elements $w,w'\in W$ by an arrow $w'\overset{s_{\alpha}}\rw w$, if $w=s_{\alpha}w'$ and $l(w)=l(w')+1$.
\end{itemize}

Let $\qP=\qL_{\qP}\N_{\qP}$ be a parabolic subgroup of $\G$ standard with respect to $\mathbf{B}$. 
The Weyl group $W_{\qL_{\qP}}$ acts as a subgroup on $W$ and the orbit space $W_{\qL_{\qP}}\backslash W$ has a distinguished set of representatives $W^{\mathbf{P}}=\{w\in W\vert w^{-1}(\Delta(\qL_{\qP},\T))\subset\Phi^+\}$, called the set of \emph{minimal coset representatives}. As a subset of $W$ it inherits the Bruhat order and thus the structure of a Hasse diagram. This structure will now be used to determine its elements.

\subsection{The actual determination of the Hasse diagram}
 Let $\Phi(\N_{\qP})$ be the set of weights for the adjoint action of $\T$ in the Lie algebra $\mathfrak{n}_{\qP}$ of the unipotent radical of $\qP$. Now, if for $w\in W$ the set $\Phi_w=\{\alpha\in\Phi^+\vert w^{-1}(\alpha)\in -\Phi^+\}$ is introduced, then
$$w^{-1}(\Delta(\qL_{\qP},\T))\subset\Phi^+\Longleftrightarrow\Phi_w\subset\Phi(\N_{\qP}),$$
simply for the fact that $\Phi^+(\qL_{\qP},\T)\cup\Phi(\N_{\qP})=\Phi^+(\G,\T)$.
So $W^{\qP}$ may as well be characterised as
$$W^{\mathbf{P}}=\{w\in W\vert \Phi_w\subset\Phi(\N_{\qP})\}.$$
The sets $\Phi_w\subset\Phi^+$ have some nice properties, two of which are stated below.

$(a)$ The cardinality of $\Phi_w$ is given as $\lvert\Phi_w\rvert=l(w)$.

$(b)$ Let $w\in W$ and $\alpha\in\Phi^+$. Then $\alpha$ is either contained in $\Phi_w$ or $\Phi_{s_{\alpha}w}$. Moreover, $l(s_{\alpha}w)>l(w)$ 
if and only if $\alpha\not\in\Phi_w$. For a simple root $\alpha$ this specifies to: $l(s_{\alpha}w)=l(w)+1$ if $\alpha\not\in\Phi_w$ and 
$l(s_{\alpha}w)=l(w)-1$ if $\alpha\in\Phi_w$.
 
The second one will be called `\textit{back-or-forth alternative}' due to the second part about simple roots.

\begin{prop}
For every $w\in W$ there are unique elements $w_{\qP}\in W_{\qL_{\qP}}$ and $w^{\qP}\in W^{\qP}$ such that
$w=w_{\qP}w^{\qP}$ and $l(w)=l(w_{\qP})+l(w^{\qP})$.
\end{prop}

\begin{rem}
The second part of the proposition exhibits $W^{\qP}$ as the set of minimal coset representatives. A proof can be found in \cite{Kos}.
\end{rem}

Let $\lambda$ be a dominant weight, dominant for the derived group of $\G$. The Bruhat order on $W$ induces a partial order on the $W$-orbit of $\lambda$: if $w'\leq w$, then  $w'(\lambda)\geq w(\lambda)$. And the idea is to realise $W^{\qP}$ as a subset of the $W$-orbit of an appropriate dominant weight $\lambda$. If $\Sigma_{\qP}$ is to denote the set of simple roots in $\Phi(\N_{\qP})$, define $\delta_{\qP}$ to be the sum of all fundamental weights corresponding to simple roots in $\Sigma_{\qP}$. The desired algorithm can be established in guise of the

\begin{prop}[Algorithm]\label{algorithm}

$(1)$ The mapping $w\mapsto w^{-1}(\delta_{\qP})$ induces a bijection between the set $W^{\qP}$ and the $W$-orbit of $\delta_{\qP}$.

$(2)$ If $w\in W^{\qP}$ and $\alpha\in\Delta$ such that $\alpha\not\in\Phi_{w^{-1}}$ and $s_{\alpha}(w^{-1}(\delta_{\qP}))\neq w^{-1}(\delta_{\qP})$,
then $w\overset{s_{w(\alpha)}}\lrw ws_{\alpha}$, $\Phi_{ws_{\alpha}}=\Phi_w\cup\{w(\alpha)\}$ and $ws_{\alpha}\in W^{\qP}$.
\end{prop}

\begin{rem}
(a) Since $s_{w(\alpha)}=ws_{\alpha}w^{-1}$, $ws_{\alpha}$ is seen to equal $s_{w(\alpha)}w$ showing the result to fit the initial definition.

(b) The first condition in part $(2)$ means not to go back, the second not to halt.
\end{rem}

\begin{rem}
A proof can be supplied with the help of \cite{B-E}.
\end{rem}

In order to actually make use of this algorithm the action of a simple reflection $s_{\alpha}$ on a weight $\lambda$ in Dynkin diagram notation has to be made explicit. 
Let $\{\alpha_1,\ldots,\alpha_k\}$ be the set of simple roots in $\Phi^+$. Let $\alpha_i^{\vee}$ denote the dual root for $\alpha_i$ and $\omega_i$ the fundamental weight corresponding to $\alpha_i$, $1\leq i\leq k$. For any weight $\lambda$ set $c_i(\lambda)=(\lambda,\alpha_i^{\vee})$ for 
the coefficient of $\lambda$ relative to $\omega_i$. In particular, $\lambda=\sum_{i=1}^kc_i(\lambda)\omega_i$ and $c_i(\lambda)$ is placed
 above the node in the Dynkin diagram corresponding to $\alpha_i$.
Now, the simple reflection $s_j=s_{\alpha_j}$ takes $\lambda$ to $s_j(\lambda)=\lambda- (\lambda,\alpha_j^{\vee})\alpha_j$, such that the coefficient
of $s_j(\lambda)$ relative to $\omega_i$ computes as
$$\begin{array}{rl}
c_i(s_j(\lambda))&=c_i(\lambda)-(\lambda,\alpha_j^{\vee})c_i(\alpha_j)\\
                                 &=c_i(\lambda)-c_j(\lambda)(\alpha_j,\alpha_i^{\vee}).
\end{array}$$
The coefficient $c_i(\lambda)$ is altered by $s_j$ if and only if $(\alpha_j,\alpha_i^{\vee})\neq0$. Suppose therefore this to be granted,
then there are two cases to be distinguished: if $j=i$, then of course $c_i(s_j(\lambda))=-c_i(\lambda)$. If on the other hand $j\neq i$, then 
$(\alpha_j,\alpha_i^{\vee})\in\{-1,-2,-3\}$ corresponding to the following pieces of Dynkin diagrams:

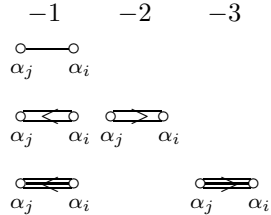
\small
$$\begin{array}{ccc}
 -1 &-2 &-3 \\
\begin{picture}(24,10)
\put(5,3){\line(1,0){16}}
\put(3,3){\makebox(0,0){$\circ$}}
\put(-1,-6){\begin{scriptsize}$\alpha_j$\end{scriptsize}}
\put(23,3){\makebox(0,0){$\circ$}}
\put(21,-6){\begin{scriptsize}$\alpha_i$\end{scriptsize}}
\end{picture} & & \\
&&\\
 \begin{picture}(24,10)
\put(4,1,4){\line(1,0){18}}
\put(4,5){\line(1,0){18}}
\put(14,3){\makebox(0,0){$<$}}
\put(3,3){\makebox(0,0){$\circ$}}
\put(-1,-6){\begin{scriptsize}$\alpha_j$\end{scriptsize}}
\put(23,3){\makebox(0,0){$\circ$}}
\put(21,-6){\begin{scriptsize}$\alpha_i$\end{scriptsize}}
\end{picture} & 
\begin{picture}(24,10)
\put(4,1,4){\line(1,0){18}}
\put(4,5){\line(1,0){18}}
\put(14,3){\makebox(0,0){$>$}}
\put(3,3){\makebox(0,0){$\circ$}}
\put(-1,-6){\begin{scriptsize}$\alpha_j$\end{scriptsize}}
\put(23,3){\makebox(0,0){$\circ$}}
\put(21,-6){\begin{scriptsize}$\alpha_i$\end{scriptsize}}
\end{picture}& \\
&&\\
\begin{picture}(24,10)
\put(4,1,4){\line(1,0){18}}
\put(5,3,2){\line(1,0){16}}
\put(4,5){\line(1,0){18}}
\put(14,3){\makebox(0,0){$<$}}
\put(3,3){\makebox(0,0){$\circ$}}
\put(-1,-6){\begin{scriptsize}$\alpha_j$\end{scriptsize}}
\put(23,3){\makebox(0,0){$\circ$}}
\put(21,-6){\begin{scriptsize}$\alpha_i$\end{scriptsize}}
\end{picture} &&
\begin{picture}(24,10)
\put(4,1,4){\line(1,0){18}}
\put(5,3,2){\line(1,0){16}}
\put(4,5){\line(1,0){18}}
\put(14,3){\makebox(0,0){$>$}}
\put(3,3){\makebox(0,0){$\circ$}}
\put(-1,-6){\begin{scriptsize}$\alpha_j$\end{scriptsize}}
\put(23,3){\makebox(0,0){$\circ$}}
\put(21,-6){\begin{scriptsize}$\alpha_i$\end{scriptsize}}
\end{picture}
\end{array}$$
\normalsize

\vspace*{2mm}

Consequently, the action of $s_j$ on a weight $\lambda$ in Dynkin diagram notation with $b=c_j(\lambda)$ is seen to be: 
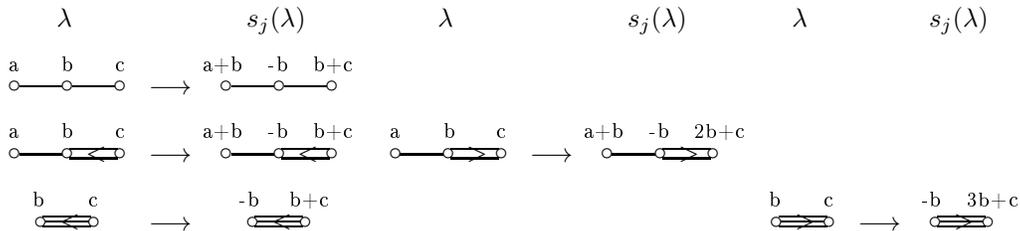
\small
$$\begin{array}{ccccccccccc}
\lambda && s_j(\lambda)&&\lambda && s_j(\lambda)&&\lambda && s_j(\lambda)\\
 &&&&&&&&&&\\
\begin{picture}(44,10)
\put(5,3){\line(1,0){16}}
\put(25,3){\line(1,0){16}}
\put(3,3){\makebox(0,0){$\circ$}}
\put(1,9){\begin{scriptsize}a\end{scriptsize}}
\put(23,3){\makebox(0,0){$\circ$}}
\put(21,9){\begin{scriptsize}b\end{scriptsize}}
\put(43,3){\makebox(0,0){$\circ$}}
\put(41,9){\begin{scriptsize}c\end{scriptsize}}
\end{picture} 
&\lrw & 
\begin{picture}(44,10)
\put(5,3){\line(1,0){16}}
\put(25,3){\line(1,0){16}}
\put(3,3){\makebox(0,0){$\circ$}}
\put(-6,9){\begin{scriptsize}a+b\end{scriptsize}}
\put(23,3){\makebox(0,0){$\circ$}}
\put(19,9){\begin{scriptsize}-b\end{scriptsize}}
\put(43,3){\makebox(0,0){$\circ$}}
\put(36,9){\begin{scriptsize}b+c\end{scriptsize}}
\end{picture}
&&&&&&&&\\
 &&&&&&&&&&\\
\begin{picture}(44,10)
\put(5,3){\line(1,0){16}}
\put(24,1,4){\line(1,0){18}}
\put(24,5){\line(1,0){18}}
\put(34,3){\makebox(0,0){$<$}}
\put(3,3){\makebox(0,0){$\circ$}}
\put(1,9){\begin{scriptsize}a\end{scriptsize}}
\put(23,3){\makebox(0,0){$\circ$}}
\put(21,9){\begin{scriptsize}b\end{scriptsize}}
\put(43,3){\makebox(0,0){$\circ$}}
\put(41,9){\begin{scriptsize}c\end{scriptsize}}
\end{picture} 
&\lrw & 
\begin{picture}(44,10)
\put(5,3){\line(1,0){16}}
\put(24,1,4){\line(1,0){18}}
\put(24,5){\line(1,0){18}}
\put(34,3){\makebox(0,0){$<$}}
\put(3,3){\makebox(0,0){$\circ$}}
\put(-6,9){\begin{scriptsize}a+b\end{scriptsize}}
\put(23,3){\makebox(0,0){$\circ$}}
\put(19,9){\begin{scriptsize}-b\end{scriptsize}}
\put(43,3){\makebox(0,0){$\circ$}}
\put(36,9){\begin{scriptsize}b+c\end{scriptsize}}
\end{picture} 
&&
\begin{picture}(44,10)
\put(5,3){\line(1,0){16}}
\put(24,1,4){\line(1,0){18}}
\put(24,5){\line(1,0){18}}
\put(34,3){\makebox(0,0){$>$}}
\put(3,3){\makebox(0,0){$\circ$}}
\put(1,9){\begin{scriptsize}a\end{scriptsize}}
\put(23,3){\makebox(0,0){$\circ$}}
\put(21,9){\begin{scriptsize}b\end{scriptsize}}
\put(43,3){\makebox(0,0){$\circ$}}
\put(41,9){\begin{scriptsize}c\end{scriptsize}}
\end{picture} 
&\lrw & 
\begin{picture}(44,10)
\put(5,3){\line(1,0){16}}
\put(24,1,4){\line(1,0){18}}
\put(24,5){\line(1,0){18}}
\put(34,3){\makebox(0,0){$>$}}
\put(3,3){\makebox(0,0){$\circ$}}
\put(-6,9){\begin{scriptsize}a+b\end{scriptsize}}
\put(23,3){\makebox(0,0){$\circ$}}
\put(19,9){\begin{scriptsize}-b\end{scriptsize}}
\put(43,3){\makebox(0,0){$\circ$}}
\put(36,9){\begin{scriptsize}2b+c\end{scriptsize}}
\end{picture} 
&&&&\\
 &&&&&&&&&&\\
\begin{picture}(24,10)
\put(4,1,4){\line(1,0){18}}
\put(5,3,2){\line(1,0){16}}
\put(4,5){\line(1,0){18}}
\put(14,3){\makebox(0,0){$<$}}
\put(3,3){\makebox(0,0){$\circ$}}
\put(0,9){\begin{scriptsize}b\end{scriptsize}}
\put(23,3){\makebox(0,0){$\circ$}}
\put(21,9){\begin{scriptsize}c\end{scriptsize}}
\end{picture} &\lrw&
\begin{picture}(24,10)
\put(4,1,4){\line(1,0){18}}
\put(5,3,2){\line(1,0){16}}
\put(4,5){\line(1,0){18}}
\put(14,3){\makebox(0,0){$<$}}
\put(3,3){\makebox(0,0){$\circ$}}
\put(-2,9){\begin{scriptsize}-b\end{scriptsize}}
\put(23,3){\makebox(0,0){$\circ$}}
\put(17,9){\begin{scriptsize}b+c\end{scriptsize}}
\end{picture} 
&&&&&&
\begin{picture}(24,10)
\put(4,1,4){\line(1,0){18}}
\put(5,3,2){\line(1,0){16}}
\put(4,5){\line(1,0){18}}
\put(14,3){\makebox(0,0){$>$}}
\put(3,3){\makebox(0,0){$\circ$}}
\put(0,9){\begin{scriptsize}b\end{scriptsize}}
\put(23,3){\makebox(0,0){$\circ$}}
\put(21,9){\begin{scriptsize}c\end{scriptsize}}
\end{picture} &\lrw&
\begin{picture}(24,10)
\put(4,1,4){\line(1,0){18}}
\put(5,3,2){\line(1,0){16}}
\put(4,5){\line(1,0){18}}
\put(14,3){\makebox(0,0){$>$}}
\put(3,3){\makebox(0,0){$\circ$}}
\put(-2,9){\begin{scriptsize}-b\end{scriptsize}}
\put(23,3){\makebox(0,0){$\circ$}}
\put(15,9){\begin{scriptsize}3b+c\end{scriptsize}}
\end{picture}
\end{array}
$$
\normalsize

The operating device for applying the algorithm may now be stated as the two stage procedure:

\textbf{Stage 1} \quad depict the crossed Dynkin diagram for $\qP$ and write the coordinates of the weight $\delta_{\qP}$ relative to the fundamental weights above the corresponding crossed, resp. uncrossed nodes. Then take all simple reflections $s_{\alpha}$ with $s_{\alpha}(\delta_{\qP})\neq\delta_{\qP}$, 
apply every one of them separately to $\delta_{\qP}$ using the rules compiled in Dynkin diagram notation above, depict the resulting weights $s_{\alpha}(\delta_{\qP})$
to the right of $\delta_{\qP}$ each joining to the latter by an arrow pointing to its transform under $s_{\alpha}$ and labelled by $s_\alpha$. 
In the next step do the analogous thing for every $s_{\alpha}(\delta_{\qP})$, that is having chosen one take all simple reflections $s_{\beta}$ with 
$\beta\not\in\Phi_{s_{\alpha}}$ and $s_{\beta}s_{\alpha}(\delta_{\qP})\neq s_{\alpha}(\delta_{\qP})$ etc. Repeat this step by step operation until
it stops, that is until it produces the situation that for any of the weights arrived at in the previous step the application of any simple reflection would 
map it to a preceding one or to itself. In terms of part $(2)$ of the proposition this situation is the one, where for the first time no simple reflection is
left meeting both hypotheses.

\textbf{Stage 2} \quad given the whole $W$-orbit of $\delta_{\qP}$ the elements of $W^{\qP}$ are obtained as the words in the simple reflections that can be built by reading the labels from the arrow to the very left to any 
arrow right to it. 

\subsubsection*{}
Finally, there shall be given two examples of Hasse diagrams. Let $\qP_2$ be the second standard maximal parabolic $\q$-subgroup of the group $\G$ from section \ref{q-form}. On the following two 
pages planar graphs are shown depicting (parts of)\footnote{The graph for $n=6$ falls short of being the actual Hasse diagram in that it lacks some arrows. Thus, it does not properly depict the partial order
on $W^{\qP}$. A refinement of the algorithm from Proposition \ref{algorithm} providing means to obtain these missing arrows will eventually be available from the book A. \v Cap, J. Slov\'ak, \emph{Parabolic geometries}, 
which is in preparation to date.} the Hasse diagram $W^{\qP_2}$ in cases $n=5,6$.

\newpage

\begin{landscape}
\noindent $\mathbf{n=5}$:

\vspace*{3mm}

$W=\langle s_1,s_2,s_3\rangle$, $\dim\mathfrak{n}_2=2\cdot5-3=7$, $\lvert W^{\qP_2}\rvert=3\cdot4=12$

\vspace*{22mm}

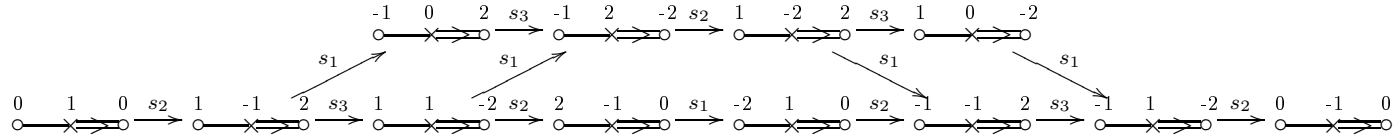
$$\def\objectstyle{\scriptstyle}
\def\labelstyle{\scriptstyle}
\xymatrix @!0 @R=12mm @C=24mm{
&&
\begin{picture}(44,10)
\put(5,3){\line(1,0){17}}
\put(25,2){\line(1,0){16}}
\put(25,4,5){\line(1,0){16}}
\put(34,3){\makebox(0,0){$>$}}
\put(3,3){\makebox(0,0){$\circ$}}
\put(1,9){\begin{scriptsize}-1\end{scriptsize}}
\put(23,3){\makebox(0,0){$\times$}}
\put(20,9){\begin{scriptsize}0\end{scriptsize}}
\put(43,3){\makebox(0,0){$\circ$}}
\put(41,9){\begin{scriptsize}2\end{scriptsize}}
\end{picture}
\ar@{->}[r]^{s_3} &
\begin{picture}(44,10)
\put(5,3){\line(1,0){17}}
\put(25,2){\line(1,0){16}}
\put(25,4,5){\line(1,0){16}}
\put(34,3){\makebox(0,0){$>$}}
\put(3,3){\makebox(0,0){$\circ$}}
\put(1,9){\begin{scriptsize}-1\end{scriptsize}}
\put(23,3){\makebox(0,0){$\times$}}
\put(20,9){\begin{scriptsize}2\end{scriptsize}}
\put(43,3){\makebox(0,0){$\circ$}}
\put(41,9){\begin{scriptsize}-2\end{scriptsize}}
\end{picture}
\ar@{->}[r]^{s_2} &
\begin{picture}(44,10)
\put(5,3){\line(1,0){17}}
\put(25,2){\line(1,0){16}}
\put(25,4,5){\line(1,0){16}}
\put(34,3){\makebox(0,0){$>$}}
\put(3,3){\makebox(0,0){$\circ$}}
\put(1,9){\begin{scriptsize}1\end{scriptsize}}
\put(23,3){\makebox(0,0){$\times$}}
\put(20,9){\begin{scriptsize}-2\end{scriptsize}}
\put(43,3){\makebox(0,0){$\circ$}}
\put(41,9){\begin{scriptsize}2\end{scriptsize}}
\end{picture}
\ar@{->}[dr]^{s_1} \ar@{->}[r]^{s_3} &
\begin{picture}(44,10)
\put(5,3){\line(1,0){17}}
\put(25,2){\line(1,0){16}}
\put(25,4,5){\line(1,0){16}}
\put(34,3){\makebox(0,0){$>$}}
\put(3,3){\makebox(0,0){$\circ$}}
\put(1,9){\begin{scriptsize}1\end{scriptsize}}
\put(23,3){\makebox(0,0){$\times$}}
\put(20,9){\begin{scriptsize}0\end{scriptsize}}
\put(43,3){\makebox(0,0){$\circ$}}
\put(41,9){\begin{scriptsize}-2\end{scriptsize}}
\end{picture}
\ar@{->}[dr]^{s_1} 
\\
\begin{picture}(44,10)
\put(5,3){\line(1,0){17}}
\put(25,2){\line(1,0){16}}
\put(25,4,5){\line(1,0){16}}
\put(34,3){\makebox(0,0){$>$}}
\put(3,3){\makebox(0,0){$\circ$}}
\put(1,9){\begin{scriptsize}0\end{scriptsize}}
\put(23,3){\makebox(0,0){$\times$}}
\put(21,9){\begin{scriptsize}1\end{scriptsize}}
\put(43,3){\makebox(0,0){$\circ$}}
\put(41,9){\begin{scriptsize}0\end{scriptsize}}
\end{picture} 
\ar@{->}[r]^{s_2} & 
\begin{picture}(44,10)
\put(5,3){\line(1,0){17}}
\put(25,2){\line(1,0){16}}
\put(25,4,5){\line(1,0){16}}
\put(34,3){\makebox(0,0){$>$}}
\put(3,3){\makebox(0,0){$\circ$}}
\put(1,9){\begin{scriptsize}1\end{scriptsize}}
\put(23,3){\makebox(0,0){$\times$}}
\put(20,9){\begin{scriptsize}-1\end{scriptsize}}
\put(43,3){\makebox(0,0){$\circ$}}
\put(41,9){\begin{scriptsize}2\end{scriptsize}}
\end{picture}
\ar@{->}[ur]^{s_1} \ar@{->}[r]^{s_3} &
\begin{picture}(44,10)
\put(5,3){\line(1,0){17}}
\put(25,2){\line(1,0){16}}
\put(25,4,5){\line(1,0){16}}
\put(34,3){\makebox(0,0){$>$}}
\put(3,3){\makebox(0,0){$\circ$}}
\put(1,9){\begin{scriptsize}1\end{scriptsize}}
\put(23,3){\makebox(0,0){$\times$}}
\put(20,9){\begin{scriptsize}1\end{scriptsize}}
\put(43,3){\makebox(0,0){$\circ$}}
\put(41,9){\begin{scriptsize}-2\end{scriptsize}}
\end{picture}
\ar@{->}[ur]^{s_1}\ar@{->}[r]^{s_2} &
\begin{picture}(44,10)
\put(5,3){\line(1,0){17}}
\put(25,2){\line(1,0){16}}
\put(25,4,5){\line(1,0){16}}
\put(34,3){\makebox(0,0){$>$}}
\put(3,3){\makebox(0,0){$\circ$}}
\put(1,9){\begin{scriptsize}2\end{scriptsize}}
\put(23,3){\makebox(0,0){$\times$}}
\put(20,9){\begin{scriptsize}-1\end{scriptsize}}
\put(43,3){\makebox(0,0){$\circ$}}
\put(41,9){\begin{scriptsize}0\end{scriptsize}}
\end{picture}
\ar@{->}[r]^{s_1} &
\begin{picture}(44,10)
\put(5,3){\line(1,0){17}}
\put(25,2){\line(1,0){16}}
\put(25,4,5){\line(1,0){16}}
\put(34,3){\makebox(0,0){$>$}}
\put(3,3){\makebox(0,0){$\circ$}}
\put(1,9){\begin{scriptsize}-2\end{scriptsize}}
\put(23,3){\makebox(0,0){$\times$}}
\put(20,9){\begin{scriptsize}1\end{scriptsize}}
\put(43,3){\makebox(0,0){$\circ$}}
\put(41,9){\begin{scriptsize}0\end{scriptsize}}
\end{picture}
\ar@{->}[r]^{s_2} &
\begin{picture}(44,10)
\put(5,3){\line(1,0){17}}
\put(25,2){\line(1,0){16}}
\put(25,4,5){\line(1,0){16}}
\put(34,3){\makebox(0,0){$>$}}
\put(3,3){\makebox(0,0){$\circ$}}
\put(1,9){\begin{scriptsize}-1\end{scriptsize}}
\put(23,3){\makebox(0,0){$\times$}}
\put(20,9){\begin{scriptsize}-1\end{scriptsize}}
\put(43,3){\makebox(0,0){$\circ$}}
\put(41,9){\begin{scriptsize}2\end{scriptsize}}
\end{picture}
\ar@{->}[r]^{s_3} &
\begin{picture}(44,10)
\put(5,3){\line(1,0){17}}
\put(25,2){\line(1,0){16}}
\put(25,4,5){\line(1,0){16}}
\put(34,3){\makebox(0,0){$>$}}
\put(3,3){\makebox(0,0){$\circ$}}
\put(1,9){\begin{scriptsize}-1\end{scriptsize}}
\put(23,3){\makebox(0,0){$\times$}}
\put(20,9){\begin{scriptsize}1\end{scriptsize}}
\put(43,3){\makebox(0,0){$\circ$}}
\put(41,9){\begin{scriptsize}-2\end{scriptsize}}
\end{picture}
\ar@{->}[r]^{s_2} &
\begin{picture}(44,10)
\put(5,3){\line(1,0){17}}
\put(25,2){\line(1,0){16}}
\put(25,4,5){\line(1,0){16}}
\put(34,3){\makebox(0,0){$>$}}
\put(3,3){\makebox(0,0){$\circ$}}
\put(1,9){\begin{scriptsize}0\end{scriptsize}}
\put(23,3){\makebox(0,0){$\times$}}
\put(20,9){\begin{scriptsize}-1\end{scriptsize}}
\put(43,3){\makebox(0,0){$\circ$}}
\put(41,9){\begin{scriptsize}0\end{scriptsize}}
\end{picture}
}
$$
\end{landscape}

\newpage
\begin{landscape}
\noindent $\mathbf{n=6}$:

\vspace*{3mm}

$W=\langle s_1,s_2,s_3,s_4\rangle$, $\dim\mathfrak{n}_2=2\cdot6-3=9$, $\lvert W^{\qP_2}\rvert=3\cdot8=24$

\begin{rem}
Notice that $s_3$ and $s_4$ commute, whereas $s_2$ and $s_4$ do not.
\end{rem}

\vspace*{6mm}

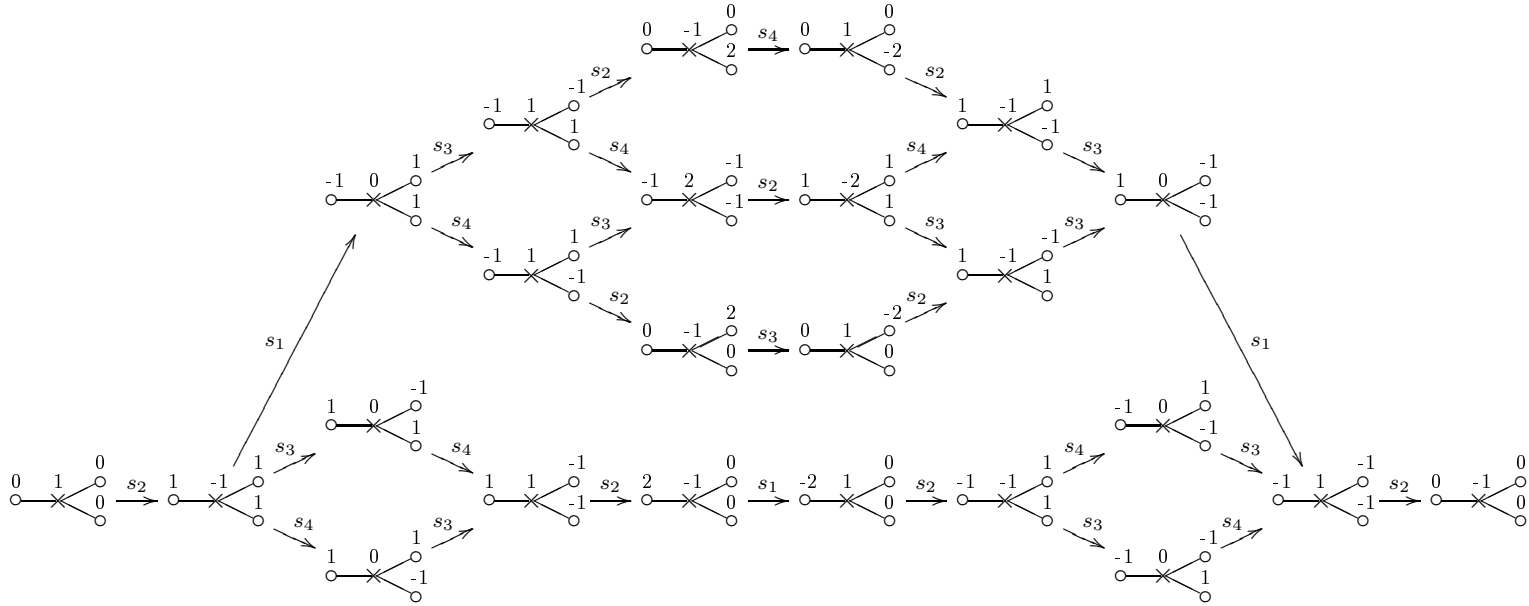
$$\def\objectstyle{\scriptstyle}
\def\labelstyle{\scriptstyle}
\xymatrix @!0 @R=10mm @C=21mm{
&&&&
\begin{picture}(38,20)
\put(5,10){\line(1,0){13}}
\put(20,10){\line(2,1){13}}
\put(20,10){\line(2,-1){13}}
\put(3,10){\makebox(0,0){$\circ$}}
\put(1,15){\begin{scriptsize}0\end{scriptsize}}
\put(19,10){\makebox(0,0){$\times$}}
\put(17,15){\begin{scriptsize}-1\end{scriptsize}}
\put(35,2){\makebox(0,0){$\circ$}}
\put(33,7){\begin{scriptsize}2\end{scriptsize}}
\put(35,17){\makebox(0,0){$\circ$}}
\put(33,22){\begin{scriptsize}0\end{scriptsize}}
\end{picture}
\ar@{->}[r]^{s_4} &
\begin{picture}(38,20)
\put(5,10){\line(1,0){13}}
\put(20,10){\line(2,1){13}}
\put(20,10){\line(2,-1){13}}
\put(3,10){\makebox(0,0){$\circ$}}
\put(1,15){\begin{scriptsize}0\end{scriptsize}}
\put(19,10){\makebox(0,0){$\times$}}
\put(17,15){\begin{scriptsize}1\end{scriptsize}}
\put(35,2){\makebox(0,0){$\circ$}}
\put(33,7){\begin{scriptsize}-2\end{scriptsize}}
\put(35,17){\makebox(0,0){$\circ$}}
\put(33,22){\begin{scriptsize}0\end{scriptsize}}
\end{picture}
\ar@{->}[dr]^{s_2}
&&&&
\\
&&&
\begin{picture}(38,20)
\put(5,10){\line(1,0){13}}
\put(20,10){\line(2,1){13}}
\put(20,10){\line(2,-1){13}}
\put(3,10){\makebox(0,0){$\circ$}}
\put(1,15){\begin{scriptsize}-1\end{scriptsize}}
\put(19,10){\makebox(0,0){$\times$}}
\put(17,15){\begin{scriptsize}1\end{scriptsize}}
\put(35,2){\makebox(0,0){$\circ$}}
\put(33,7){\begin{scriptsize}1\end{scriptsize}}
\put(35,17){\makebox(0,0){$\circ$}}
\put(33,22){\begin{scriptsize}-1\end{scriptsize}}
\end{picture}
\ar@{->}[ur]^{s_2} 
\ar@{->}[dr]^{s_4} &&&
\begin{picture}(38,20)
\put(5,10){\line(1,0){13}}
\put(20,10){\line(2,1){13}}
\put(20,10){\line(2,-1){13}}
\put(3,10){\makebox(0,0){$\circ$}}
\put(1,15){\begin{scriptsize}1\end{scriptsize}}
\put(19,10){\makebox(0,0){$\times$}}
\put(17,15){\begin{scriptsize}-1\end{scriptsize}}
\put(35,2){\makebox(0,0){$\circ$}}
\put(33,7){\begin{scriptsize}-1\end{scriptsize}}
\put(35,17){\makebox(0,0){$\circ$}}
\put(33,22){\begin{scriptsize}1\end{scriptsize}}
\end{picture}
\ar@{->}[dr]^{s_3}
&&&
\\
&&
\begin{picture}(38,20)
\put(5,10){\line(1,0){13}}
\put(20,10){\line(2,1){13}}
\put(20,10){\line(2,-1){13}}
\put(3,10){\makebox(0,0){$\circ$}}
\put(1,15){\begin{scriptsize}-1\end{scriptsize}}
\put(19,10){\makebox(0,0){$\times$}}
\put(17,15){\begin{scriptsize}0\end{scriptsize}}
\put(35,2){\makebox(0,0){$\circ$}}
\put(33,7){\begin{scriptsize}1\end{scriptsize}}
\put(35,17){\makebox(0,0){$\circ$}}
\put(33,22){\begin{scriptsize}1\end{scriptsize}}
\end{picture}
\ar@{->}[ur]^{s_3} 
\ar@{->}[dr]^{s_4} &&
\begin{picture}(38,20)
\put(5,10){\line(1,0){13}}
\put(20,10){\line(2,1){13}}
\put(20,10){\line(2,-1){13}}
\put(3,10){\makebox(0,0){$\circ$}}
\put(1,15){\begin{scriptsize}-1\end{scriptsize}}
\put(19,10){\makebox(0,0){$\times$}}
\put(17,15){\begin{scriptsize}2\end{scriptsize}}
\put(35,2){\makebox(0,0){$\circ$}}
\put(33,7){\begin{scriptsize}-1\end{scriptsize}}
\put(35,17){\makebox(0,0){$\circ$}}
\put(33,22){\begin{scriptsize}-1\end{scriptsize}}
\end{picture}
\ar@{->}[r]^{s_2} &
\begin{picture}(38,20)
\put(5,10){\line(1,0){13}}
\put(20,10){\line(2,1){13}}
\put(20,10){\line(2,-1){13}}
\put(3,10){\makebox(0,0){$\circ$}}
\put(1,15){\begin{scriptsize}1\end{scriptsize}}
\put(19,10){\makebox(0,0){$\times$}}
\put(17,15){\begin{scriptsize}-2\end{scriptsize}}
\put(35,2){\makebox(0,0){$\circ$}}
\put(33,7){\begin{scriptsize}1\end{scriptsize}}
\put(35,17){\makebox(0,0){$\circ$}}
\put(33,22){\begin{scriptsize}1\end{scriptsize}}
\end{picture}
\ar@{->}[ur]^{s_4} 
\ar@{->}[dr]^{s_3} &&
\begin{picture}(38,20)
\put(5,10){\line(1,0){13}}
\put(20,10){\line(2,1){13}}
\put(20,10){\line(2,-1){13}}
\put(3,10){\makebox(0,0){$\circ$}}
\put(1,15){\begin{scriptsize}1\end{scriptsize}}
\put(19,10){\makebox(0,0){$\times$}}
\put(17,15){\begin{scriptsize}0\end{scriptsize}}
\put(35,2){\makebox(0,0){$\circ$}}
\put(33,7){\begin{scriptsize}-1\end{scriptsize}}
\put(35,17){\makebox(0,0){$\circ$}}
\put(33,22){\begin{scriptsize}-1\end{scriptsize}}
\end{picture}
\ar@{->}[ddddr]^{s_1} 
&&
\\
&&&
\begin{picture}(38,20)
\put(5,10){\line(1,0){13}}
\put(20,10){\line(2,1){13}}
\put(20,10){\line(2,-1){13}}
\put(3,10){\makebox(0,0){$\circ$}}
\put(1,15){\begin{scriptsize}-1\end{scriptsize}}
\put(19,10){\makebox(0,0){$\times$}}
\put(17,15){\begin{scriptsize}1\end{scriptsize}}
\put(35,2){\makebox(0,0){$\circ$}}
\put(33,7){\begin{scriptsize}-1\end{scriptsize}}
\put(35,17){\makebox(0,0){$\circ$}}
\put(33,22){\begin{scriptsize}1\end{scriptsize}}
\end{picture}
\ar@{->}[ur]^{s_3} 
\ar@{->}[dr]^{s_2} &&&
\begin{picture}(38,20)
\put(5,10){\line(1,0){13}}
\put(20,10){\line(2,1){13}}
\put(20,10){\line(2,-1){13}}
\put(3,10){\makebox(0,0){$\circ$}}
\put(1,15){\begin{scriptsize}1\end{scriptsize}}
\put(19,10){\makebox(0,0){$\times$}}
\put(17,15){\begin{scriptsize}-1\end{scriptsize}}
\put(35,2){\makebox(0,0){$\circ$}}
\put(33,7){\begin{scriptsize}1\end{scriptsize}}
\put(35,17){\makebox(0,0){$\circ$}}
\put(33,22){\begin{scriptsize}-1\end{scriptsize}}
\end{picture}
\ar@{->}[ur]^{s_3}
&&&
\\
&&&&
\begin{picture}(38,20)
\put(5,10){\line(1,0){13}}
\put(20,10){\line(2,1){13}}
\put(20,10){\line(2,-1){13}}
\put(3,10){\makebox(0,0){$\circ$}}
\put(1,15){\begin{scriptsize}0\end{scriptsize}}
\put(19,10){\makebox(0,0){$\times$}}
\put(17,15){\begin{scriptsize}-1\end{scriptsize}}
\put(35,2){\makebox(0,0){$\circ$}}
\put(33,7){\begin{scriptsize}0\end{scriptsize}}
\put(35,17){\makebox(0,0){$\circ$}}
\put(33,22){\begin{scriptsize}2\end{scriptsize}}
\end{picture}
\ar@{->}[r]^{s_3} &
\begin{picture}(38,20)
\put(5,10){\line(1,0){13}}
\put(20,10){\line(2,1){13}}
\put(20,10){\line(2,-1){13}}
\put(3,10){\makebox(0,0){$\circ$}}
\put(1,15){\begin{scriptsize}0\end{scriptsize}}
\put(19,10){\makebox(0,0){$\times$}}
\put(17,15){\begin{scriptsize}1\end{scriptsize}}
\put(35,2){\makebox(0,0){$\circ$}}
\put(33,7){\begin{scriptsize}0\end{scriptsize}}
\put(35,17){\makebox(0,0){$\circ$}}
\put(33,22){\begin{scriptsize}-2\end{scriptsize}}
\end{picture}
\ar@{->}[ur]^{s_2} 
&&&&
\\
&&
\begin{picture}(38,20)
\put(5,10){\line(1,0){13}}
\put(20,10){\line(2,1){13}}
\put(20,10){\line(2,-1){13}}
\put(3,10){\makebox(0,0){$\circ$}}
\put(1,15){\begin{scriptsize}1\end{scriptsize}}
\put(19,10){\makebox(0,0){$\times$}}
\put(17,15){\begin{scriptsize}0\end{scriptsize}}
\put(35,2){\makebox(0,0){$\circ$}}
\put(33,7){\begin{scriptsize}1\end{scriptsize}}
\put(35,17){\makebox(0,0){$\circ$}}
\put(33,22){\begin{scriptsize}-1\end{scriptsize}}
\end{picture}
\ar@{->}[dr]^{s_4}
&&&&&
\begin{picture}(38,20)
\put(5,10){\line(1,0){13}}
\put(20,10){\line(2,1){13}}
\put(20,10){\line(2,-1){13}}
\put(3,10){\makebox(0,0){$\circ$}}
\put(1,15){\begin{scriptsize}-1\end{scriptsize}}
\put(19,10){\makebox(0,0){$\times$}}
\put(17,15){\begin{scriptsize}0\end{scriptsize}}
\put(35,2){\makebox(0,0){$\circ$}}
\put(33,7){\begin{scriptsize}-1\end{scriptsize}}
\put(35,17){\makebox(0,0){$\circ$}}
\put(33,22){\begin{scriptsize}1\end{scriptsize}}
\end{picture}
\ar@{->}[dr]^{s_3} &&
\\
\begin{picture}(38,20)
\put(5,10){\line(1,0){13}}
\put(20,10){\line(2,1){13}}
\put(20,10){\line(2,-1){13}}
\put(3,10){\makebox(0,0){$\circ$}}
\put(1,15){\begin{scriptsize}0\end{scriptsize}}
\put(19,10){\makebox(0,0){$\times$}}
\put(17,15){\begin{scriptsize}1\end{scriptsize}}
\put(35,2){\makebox(0,0){$\circ$}}
\put(33,7){\begin{scriptsize}0\end{scriptsize}}
\put(35,17){\makebox(0,0){$\circ$}}
\put(33,22){\begin{scriptsize}0\end{scriptsize}}
\end{picture}
\ar@{->}[r]^{s_2}  &
\begin{picture}(38,20)
\put(5,10){\line(1,0){13}}
\put(20,10){\line(2,1){13}}
\put(20,10){\line(2,-1){13}}
\put(3,10){\makebox(0,0){$\circ$}}
\put(1,15){\begin{scriptsize}1\end{scriptsize}}
\put(19,10){\makebox(0,0){$\times$}}
\put(17,15){\begin{scriptsize}-1\end{scriptsize}}
\put(35,2){\makebox(0,0){$\circ$}}
\put(33,7){\begin{scriptsize}1\end{scriptsize}}
\put(35,17){\makebox(0,0){$\circ$}}
\put(33,22){\begin{scriptsize}1\end{scriptsize}}
\end{picture}
\ar@{->}[uuuur]^{s_1} 
\ar@{->}[ur]^{s_3} 
\ar@{->}[dr]^{s_4} &&
\begin{picture}(38,20)
\put(5,10){\line(1,0){13}}
\put(20,10){\line(2,1){13}}
\put(20,10){\line(2,-1){13}}
\put(3,10){\makebox(0,0){$\circ$}}
\put(1,15){\begin{scriptsize}1\end{scriptsize}}
\put(19,10){\makebox(0,0){$\times$}}
\put(17,15){\begin{scriptsize}1\end{scriptsize}}
\put(35,2){\makebox(0,0){$\circ$}}
\put(33,7){\begin{scriptsize}-1\end{scriptsize}}
\put(35,17){\makebox(0,0){$\circ$}}
\put(33,22){\begin{scriptsize}-1\end{scriptsize}}
\end{picture}
\ar@{->}[r]^{s_2}  &
\begin{picture}(38,20)
\put(5,10){\line(1,0){13}}
\put(20,10){\line(2,1){13}}
\put(20,10){\line(2,-1){13}}
\put(3,10){\makebox(0,0){$\circ$}}
\put(1,15){\begin{scriptsize}2\end{scriptsize}}
\put(19,10){\makebox(0,0){$\times$}}
\put(17,15){\begin{scriptsize}-1\end{scriptsize}}
\put(35,2){\makebox(0,0){$\circ$}}
\put(33,7){\begin{scriptsize}0\end{scriptsize}}
\put(35,17){\makebox(0,0){$\circ$}}
\put(33,22){\begin{scriptsize}0\end{scriptsize}}
\end{picture}
\ar@{->}[r]^{s_1}  &
\begin{picture}(38,20)
\put(5,10){\line(1,0){13}}
\put(20,10){\line(2,1){13}}
\put(20,10){\line(2,-1){13}}
\put(3,10){\makebox(0,0){$\circ$}}
\put(1,15){\begin{scriptsize}-2\end{scriptsize}}
\put(19,10){\makebox(0,0){$\times$}}
\put(17,15){\begin{scriptsize}1\end{scriptsize}}
\put(35,2){\makebox(0,0){$\circ$}}
\put(33,7){\begin{scriptsize}0\end{scriptsize}}
\put(35,17){\makebox(0,0){$\circ$}}
\put(33,22){\begin{scriptsize}0\end{scriptsize}}
\end{picture}
\ar@{->}[r]^{s_2}  &
\begin{picture}(38,20)
\put(5,10){\line(1,0){13}}
\put(20,10){\line(2,1){13}}
\put(20,10){\line(2,-1){13}}
\put(3,10){\makebox(0,0){$\circ$}}
\put(1,15){\begin{scriptsize}-1\end{scriptsize}}
\put(19,10){\makebox(0,0){$\times$}}
\put(17,15){\begin{scriptsize}-1\end{scriptsize}}
\put(35,2){\makebox(0,0){$\circ$}}
\put(33,7){\begin{scriptsize}1\end{scriptsize}}
\put(35,17){\makebox(0,0){$\circ$}}
\put(33,22){\begin{scriptsize}1\end{scriptsize}}
\end{picture}
\ar@{->}[ur]^{s_4} 
\ar@{->}[dr]^{s_3} &&
\begin{picture}(38,20)
\put(5,10){\line(1,0){13}}
\put(20,10){\line(2,1){13}}
\put(20,10){\line(2,-1){13}}
\put(3,10){\makebox(0,0){$\circ$}}
\put(1,15){\begin{scriptsize}-1\end{scriptsize}}
\put(19,10){\makebox(0,0){$\times$}}
\put(17,15){\begin{scriptsize}1\end{scriptsize}}
\put(35,2){\makebox(0,0){$\circ$}}
\put(33,7){\begin{scriptsize}-1\end{scriptsize}}
\put(35,17){\makebox(0,0){$\circ$}}
\put(33,22){\begin{scriptsize}-1\end{scriptsize}}
\end{picture}
\ar@{->}[r]^{s_2} &
\begin{picture}(38,20)
\put(5,10){\line(1,0){13}}
\put(20,10){\line(2,1){13}}
\put(20,10){\line(2,-1){13}}
\put(3,10){\makebox(0,0){$\circ$}}
\put(1,15){\begin{scriptsize}0\end{scriptsize}}
\put(19,10){\makebox(0,0){$\times$}}
\put(17,15){\begin{scriptsize}-1\end{scriptsize}}
\put(35,2){\makebox(0,0){$\circ$}}
\put(33,7){\begin{scriptsize}0\end{scriptsize}}
\put(35,17){\makebox(0,0){$\circ$}}
\put(33,22){\begin{scriptsize}0\end{scriptsize}}
\end{picture}
\\
&&
\begin{picture}(38,20)
\put(5,10){\line(1,0){13}}
\put(20,10){\line(2,1){13}}
\put(20,10){\line(2,-1){13}}
\put(3,10){\makebox(0,0){$\circ$}}
\put(1,15){\begin{scriptsize}1\end{scriptsize}}
\put(19,10){\makebox(0,0){$\times$}}
\put(17,15){\begin{scriptsize}0\end{scriptsize}}
\put(35,2){\makebox(0,0){$\circ$}}
\put(33,7){\begin{scriptsize}-1\end{scriptsize}}
\put(35,17){\makebox(0,0){$\circ$}}
\put(33,22){\begin{scriptsize}1\end{scriptsize}}
\end{picture}
\ar@{->}[ur]^{s_3} 
&&&&&
\begin{picture}(38,20)
\put(5,10){\line(1,0){13}}
\put(20,10){\line(2,1){13}}
\put(20,10){\line(2,-1){13}}
\put(3,10){\makebox(0,0){$\circ$}}
\put(1,15){\begin{scriptsize}-1\end{scriptsize}}
\put(19,10){\makebox(0,0){$\times$}}
\put(17,15){\begin{scriptsize}0\end{scriptsize}}
\put(35,2){\makebox(0,0){$\circ$}}
\put(33,7){\begin{scriptsize}1\end{scriptsize}}
\put(35,17){\makebox(0,0){$\circ$}}
\put(33,22){\begin{scriptsize}-1\end{scriptsize}}
\end{picture}
\ar@{->}[ur]^{s_4} &
}
$$
\end{landscape}

\end{appendix}


\begin{thebibliography}{99}
\bibitem[B-E]{B-E}R. Baston, M. Eastwood, \emph{The Penrose transform. Its interaction with representation theory}. Oxford Mathematical Monographs. Oxford Science Publications. The Clarendon Press, Oxford University Press, New York, 1989. xvi+213 pp.
\bibitem[B-W]{B-W}A. Borel, N. Wallach, \emph{Continuous cohomology, discrete subgroups and representations of reductive groups}. Second edition. Mathematical Surveys and Monographs, 67. American Mathematical Society, Providence, RI, 2000. xviii+260 pp.
\bibitem[Bou]{Bou}N. Bourbaki, \emph{Groupes et alg\`ebres de Lie}. Cha IV-VI. Hermann, Paris, 1968.
\bibitem[Fra]{Fra}J. Franke, \emph{Harmonic analysis in weighted $L\sb 2$-spaces}. Ann. Sci. \'Ecole Norm. Sup. (4) 
\textbf{31} (1998),  no. 2, 181--279.
\bibitem[F-S]{F-S}J. Franke, J. Schwermer, \emph{A decomposition of spaces of automorphic forms and the Eisenstein cohomology of arithmetic groups}. Math. Ann. \textbf{311} (1998), No.4, 765--790.
\bibitem[Har1]{Har1}G. Harder, \emph{On the cohomology of $SL_2(\mathcal{O})$}.  Lie groups and their representations (Proc. Summer School on Group Representations of the Bolyai J\'anos Math. Soc., Budapest, 1971),  pp. 139--150. Halsted, New York, 1975.
\bibitem[Har2]{Har2}G. Harder, \emph{On the cohomology of discrete arithmetically defined groups}.
Discrete Subgroups of Lie Groups Appl. Moduli, Pap. Bombay Colloq. 1973 (1975), 129--160.
\bibitem[Har3]{Har3}G. Harder, \emph{Eisenstein cohomology of arithmetic groups: the case $GL_2$}. Invent. math. 
\textbf{89} (1987), 37--118. 
\bibitem[Ha-Ch]{Ha-Ch}Harish-Chandra, \emph{Automorphic forms on semisimple Lie groups}. Notes by J. G. M. Mars. Lecture Notes in Mathematics, No. 62 Springer-Verlag, Berlin-New York 1968. x+138 pp.
\bibitem[H-S]{Ha-Sch}T. Hayata, J. Schwermer, \emph{On arithmetic subgroups of a $\q$-rank $2$ form of $SU(2,2)$ and their automorphic cohomology}. J. Math. Soc. Japan \textbf{57} (2005), No. 2, 357--385.
\bibitem[Kos]{Kos}B. Kostant, \emph{Lie algebra cohomology and the generalized Borel-Weil theorem}.  Ann. of Math. (2) \textbf{74} (1961), 329--387.
\textbf{86} (1997),  no. 1, 39--78.
\bibitem[Lan]{Lan}R. Langlands, \emph{On the functional equations satisfied by Eisenstein series}. Lecture Notes in Mathematics, Vol. 544. Springer-Verlag, Berlin-New York, 1976. v+337 pp.
\bibitem[L-S]{L-S}J.-S. Li, J. Schwermer, \emph{On the Eisenstein cohomology of arithmetic groups}.  
Duke Math. J.  \textbf{123} (2004), no. 1, 141--169.
\bibitem[M-W]{MW}C. Moeglin, J.-L. Waldspurger, \emph{Spectral decomposition and Eisenstein series}. Une paraphrase de l'\'ecriture. Cambridge Tracts in Mathematics, 113. Cambridge University Press, Cambridge, 1995. xxviii+338 pp.
\bibitem[R-S]{RS}J. Rohlfs, B. Speh, \emph{Representations with cohomology in the discrete spectrum of subgroups of $SO(n,1)(\mathbb{Z})$ and Lefschetz numbers}.  Ann. Sci. \'Ecole Norm. Sup. (4)  20  (1987),  no. 1, 89--136.
\bibitem[Sch1]{Sch1}J. Schwermer, \emph{Kohomologie arithmetisch definierter Gruppen und Eisensteinreihen}. Lecture Notes in Mathematics, 988. Springer-Verlag, Berlin, 1983. iv+170 pp.
\bibitem[Sch2]{Sch2}J. Schwermer, \emph{On arithmetic quotients of the Siegel upper half space of degree two}. Compositio Math. \textbf{58} (1986), 233--258.
\bibitem[Sch3]{Sch3}J. Schwermer, \emph{Cohomology of arithmetic groups, automorphic forms and $L$-functions}.  Cohomology of arithmetic groups and automorphic forms (Luminy-Marseille, 1989),  1--29, Lecture Notes in Math., 1447, Springer, Berlin, 1990.
\bibitem[Sch4]{Sch4}J. Schwermer, \emph{Eisenstein series and cohomology of arithmetic groups: The generic case}. Invent. Math. \textbf{116} (1994), No.1-3, 481--511.
\bibitem[V-Z]{VZ}D. Vogan, G. Zuckerman, \emph{Unitary representations with non-zero cohomology}. Compositio Math.  \textbf{53}  (1984),  no. 1, 51--90.


\end{thebibliography}
\end{document}